%% file: ms.tex
\documentclass[letterpaper]{article}

\usepackage{arxiv}

\usepackage{amsmath, amssymb, latexsym, multicol, amsthm, bm, placeins}
\usepackage[export]{adjustbox}
\usepackage{color}
\usepackage{multirow}

\usepackage{style}

\makeatletter\@addtoreset{equation}{section}\makeatother
\newcommand{\Hr}{{\rm H}}
\newcommand{\Wr}{{\rm W}}
\newcommand{\Lr}{{\rm L}}
\newcommand{\Cr}{{\rm C}}

\newcommand{\cblue}[1]{{\leavevmode\color{black}#1}}

\title{Variational Physics Informed Neural Networks: the role of quadratures and test functions}

\author{Stefano Berrone\thanks{Dipartimento di Scienze Matematiche, Politecnico di Torino, Corso Duca degli Abruzzi 24, 10129 Torino, Italy. stefano.berrone@polito.it (S. Berrone), claudio.canuto@polito.it (C. Canuto), moreno.pintore@polito.it (M. Pintore).}
	\And
	Claudio Canuto\footnotemark[1]
	\And Moreno Pintore\footnotemark[1]
}

\renewcommand{\shorttitle}{VPINNs: the role of quadratures and test functions}

\begin{document}

\maketitle

\begin{abstract}
In this work we analyze how quadrature rules of different precisions and piecewise polynomial test functions of different degrees affect the convergence rate of Variational Physics Informed Neural Networks (VPINN) with respect to mesh refinement, while solving elliptic boundary-value problems. Using a Petrov-Galerkin framework relying on an inf-sup condition, we derive an a priori error estimate in the energy norm between the exact solution and a suitable high-order piecewise interpolant of a computed neural network. Numerical experiments confirm the theoretical predictions and highlight the importance of the inf-sup condition. Our results suggest, somehow counterintuitively, that for smooth solutions the best strategy to achieve a high decay rate of the error consists in choosing test functions of the lowest polynomial degree, while using quadrature formulas of suitably high precision.
\end{abstract}

\vspace{0.3cm}
\noindent \textbf{\textit{MSC-class}} 35B45, 35J20, 35Q93, 65K10, 65N20, 68T07

\section{Introduction}
\label{intro_section}

\input{introduction}

\input{VPINN_apriori_estimates}

\input{implementation}

\input{numerical_results}

\input{conclusions}




\bibliographystyle{siam}
\bibliography{bibliography}


\appendix
\renewcommand{\thesection}{\Alph{section}.\arabic{section}}
\input{interpolation_construction}

\end{document}

%% file: introduction.tex
Exploiting the recent advances in artificial intelligence and, in particular, in deep learning, several innovative numerical techniques have been developed in the last few years to compute numerical solutions of partial differential equations (PDEs). In such methods, the solution is approximated by a neural network that is trained by taking advantage of the knowledge of the underlying differential equation. One of the earliest models involving a neural network was described in \cite{raissi2019physics}: it is based on the concept of Physics Informed Neural Networks (PINN) and it inspired further works such as e.g. \cite{tartakovsky2018learning} or \cite{yang2019adversarial}, until the recent paper \cite{LanthalerMishraKarniadakis2021} which presents a very general framework for the solution of operator equations by deep neural networks.

In such papers, given an arbitrary PDE coupled with proper boundary conditions, the training of the PINN aims at finding the weights ${\mathbf w}$ of a neural network such that the associated function $u^{\NN}(\x;{\mathbf w})$ minimizes some functional of the equation residual while satisfying as much as possible the imposed boundary conditions. To do so, the neural network is trained to minimize the residual only at a finite set of collocation points and additional terms are added to the loss function in order to force the network to approximately satisfy the boundary conditions. Thanks to the good approximation properties of neural networks, formally proved e.g. in \cite{elbrachter2021deep}, \cite{guhring2020error}, \cite{opschoor2020deep}, \cite{kutyniok2021theoretical}, \cite{opschoor2021exponential} and \cite{gonon2021deep} under suitable assumptions, the PINN approach looks very promising because it is able to efficiently and accurately compute approximate solutions of arbitrary PDEs encoding their structures in the loss function. 

Subsequently, the PINN paradigm has been further developed in \cite{kharazmi2019variational} to obtain the so-called Variational Physics Informed Neural Networks (VPINN). The main differences with respect to the PINN are that the weak formulation of the PDE is exploited, the collocation points are replaced by test functions, and quadrature points are used to compute the integrals involved in the variational residuals. In such a method the solution is still approximated by a neural network, but the test functions are represented by a finite set of known functions or by a second neural network (see \cite{zang2020weak}); therefore, the technique can be seen as a Petrov-Galerkin method.  The method is more flexible than the standard PINN because the integration by parts, involved in the weak formulation, decreases the required regularity of the approximate solution. Furthermore, the fact that the dataset used in the training phase consists of quadrature points significantly reduces the computational cost of the training phase. Indeed, quadrature points are, in general, much fewer than collocation points. 

Combining the VPINN with the Finite Element Method (FEM), the authors of \cite{khodayi2020varnet} developed VarNet, a VPINN that exploits the test functions of the $\mathbb{P}_1$-FEM. Such a work has been then extended in \cite{kharazmi2021hp} to consider arbitrary high-order polynomials as test functions, as in the $hp$ version of the FEM.  
Although the authors of the cited works empirically observed that both PINNs and VPINNs are able to efficiently approximate the desired solution, \cblue{no proof of a priori error estimates  with convergence rates is provided for VPINNs. On the contrary, rigorous a posteriori error analyses are already available (see, for instance, \cite{Mishra2021estimates}). Recently (see \cite{BeCaPi2022aposteriori}) we derived a posteriori error estimates for the discretization setting considered in this paper.
}

The purpose of this paper is to investigate how the choice of piecewise polynomial test functions and quadrature formulas influence the accuracy of the resulting VPINN approximation of second-order elliptic problems. One might think that test functions of high polynomial degree are needed to get a high order of accuracy; we prove that this is not the case, actually we indicate that precisely the opposite is true: it is more convenient to keep the degree of the test functions as low as possible, while using quadrature formulas of precision as high as possible. Indeed, for sufficiently smooth solutions, the error decay rate is given by
$$
q+2 - k_\text{test} \,,
$$
where $q$ is the precision of the quadrature formula and $k_\text{test}$ is the degree of the test functions. 

Using a Petrov-Galerkin framework, we derive an a priori error estimate in the energy norm between the exact solution and a suitable piecewise polynomial interpolant of the computed neural network; we assume that the architecture of the neural network is fixed and sufficiently rich, and we explore the behaviour of the error versus the size of the mesh supporting the test functions. Our analysis relies upon the validity of an inf-sup condition between the spaces of test functions and the space in which the neural network is interpolated. 

Numerical experiments confirm the theoretical prediction. Interestingly, in our experiments the error between the exact solution and the computed neural network decays asymptotically with the same rate as predicted by our theory for the interpolant of the network; however, this behaviour cannot be rigorously guaranteed, since in general the minimization problem which defines the computed neural network is underdetermined, and the computed neural network may be affected by spurious components. Indeed, we show that for a problem with zero data the minimization of the loss function may yield non-vanishing neural networks. With the method proposed in this paper, we combine the efficiency of the VPINN approach with the availability of a sound and certified convergence analysis.

The paper is organized as follows. In Sect. \ref{sec:setting} we introduce the elliptic problem we are focusing on, and we also present the way in which the Dirichlet boundary conditions are exactly imposed, which is uncommon in PINNs and VPINNs but can be generalized as in \cite{sukumar2021exact}. In Sect. \ref{sec:discretization} we focus on the numerical discretization; in particular, the involved neural network architecture is described in Sect. \ref{sec:networks}, while the problem discretization and the corresponding loss function are described in Sect. \ref{sec:sub_discretization}. Here we also introduce an interpolation operator $\mathcal I_H$ applied to the neural networks.
Sect. \ref{sec:error_estimates} is the key \cblue{theoretical} section: through a series of preliminary results, we formally derive the a priori error estimate, the main result being Theorem \ref{teo:a-priori bound}.
In Sect. \ref{sec:implementation} we specify the parameters of the neural network used for the numerical tests and the training phase details. We also analyse the consequences of fulfilling the inf-sup condition in connection with the VPINN efficiency. \cblue{Various numerical tests are presented and discussed in Sect. \ref{sect:num-res} for two-dimensional elliptic problems. Such tests empirically confirm the validity of the a priori estimate in different scenarios. Furthermore,} \cblue{we compare the accuracy of the proposed method with that of a standard PINN and the non-interpolated VPINN.} \cblue{ We also analyse the relationship between the neural network hyperparameters and the VPINN accuracy, and we highlight, with numerical experiments and analytical examples, the importance of the inf-sup condition. 
In Sect. \ref{sect:parametric}, we show that our VPINN can be easily adapted to solve a parametric nonlinear PDE, with accurate results for the whole range of parameters.}
Finally, in Sect. \ref{sec:conclusions}, we draw some conclusions and highlight the future perspective of the current work.

%% file: VPINN_apriori_estimates.tex
\section{The model boundary-value problem} \label{sec:setting}

\newcommand{\kt}{k_\text{int}}

Let $\Omega \subset \mathbb{R}^n$ be  a bounded polygonal/polyhedral domain with boundary $\Gamma=\partial\Omega$, partitioned into $\Gamma=\Gamma_D \cup \Gamma_N$ with $\Gamma_D \cap \Gamma_N = \emptyset$ and $\text{meas}_{n-1}(\Gamma_D)>0$.

Let us consider the model elliptic boundary-value problem
\begin{equation}\label{eq:model-pb}
\begin{cases}
Lu:=-\nabla \cdot (\mu \nabla u) + \boldsymbol{\beta}\cdot \nabla u + \sigma u =f & \text{in \ } \Omega\,, \\
u=g & \text{on \ } \Gamma_D \,, \\
\mu \frac{\partial u}{\partial n} = \psi  & \text{on \ } \Gamma_N \,, 
\end{cases}
\end{equation}
where $\mu, \sigma \in \Lr^\infty(\Omega)$, $ \boldsymbol{\beta} \in (\Wr^{1,\infty}(\Omega))^n$ satisfy $\mu \geq \mu_0$, $\sigma - \frac12 \nabla \cdot \boldsymbol{\beta} \geq 0$ in $\Omega$ for some constant $\mu_0>0$, whereas $f \in L^2(\Omega)$, $g=\bar{u}_{\vert\Gamma_D}$ for some $\bar{u} \in \Hr^1(\Omega)$, and $\psi \in \Lr^2(\Gamma_N)$.

\smallskip
Define the spaces $U=\Hr^1(\Omega)$, $V=\Hr^1_{0,\Gamma_D}(\Omega):=\{v \in U : v_{\vert\Gamma_D}=0 \}$,  the bilinear form $a:U\times V \to \mathbb{R}$ and the linear form $F:V\to \mathbb{R}$ such that 
\begin{equation}\label{eq:forms a-F}
a(w,v)=\int_\Omega \mu \nabla w \cdot \nabla v + \boldsymbol{\beta}\cdot \nabla w \, v + \sigma w \, v\,, \qquad F(v)=\int_\Omega f \, v + \int_{\Gamma_N} \psi \, v \,;
\end{equation}
denote by $\alpha \geq \mu_0$  the coercivity constant of the form $a$, and by $\Vert a \Vert$, $\Vert F \Vert$ the continuity constants of the forms $a$ and $F$. Problem \eqref{eq:model-pb} is formulated variationally as follows: {\it Find $u \in \bar{u}+V $ such that}
\begin{equation}\label{eq:model-pb-var}
a(u,v)=F(v) \qquad \forall v \in V\,.
\end{equation}
We assume that we can represent $u$ in the form 
\begin{equation}\label{eq:split-u}
u= \bar{u}+\Phi \tilde{u} \,,
\end{equation}
for some (known) smooth function $\Phi \in V$ and some $\tilde{u} \in U$ having the same smoothness of $u$.
Let us introduce the affine mapping
\begin{equation}\label{eq:def-B}
B :U \to \bar{u}+V \qquad \text{ such that} \quad Bw=\bar{u}+\Phi w
\end{equation} 
which enforces the given Dirichlet boundary condition. Then, Problem \eqref{eq:model-pb-var} can be equivalently formulated as follows: {\it Find $\tilde{u} \in U$ such that}
\begin{equation}\label{eq:model-pb-var-2}
a(B\tilde{u},v)=F(v) \qquad \forall v \in V\,.
\end{equation}

\begin{remark}[Enforcement of the Dirichlet conditions]\label{rem:bcs}{\rm
The approach we follow to enforce Dirichlet boundary conditions will allow us to deal with a loss function which is built solely by the residuals of the variational equations. Other approaches are obviously possible: for instance, one could augment such loss function by a term penalizing the distance of the numerical solution from the data on $\Gamma_D$, or adopt a Nitsche's type variational formulation of the boundary-value problem \cite{nitscheweak}. Both strategies involve parameters which may need a tuning, whereas in our approach the definition of the loss function is simple and natural, allowing us to focus on the performances of the neural networks.

}
\end{remark}

\section{The VPINN-based numerical discretization} \label{sec:discretization}

In this section, we first introduce the class of neural networks used in this paper, then we describe the numerical discretization of the boundary-value problem \eqref{eq:model-pb-var-2}, which uses neural networks to represent the discrete solution and piecewise polynomial functions to enforce the variational equations. An inf-sup stable Petrov-Galerkin formulation is introduced which guarantees stability and convergence, as indicated in Sect. \ref{sec:error_estimates}; this is the main difference between the proposed method and other formulations, such as \cite{kharazmi2019variational, kharazmi2021hp}.


\subsection{Neural networks} \label{sec:networks}

In this work we only use fully-connected feed-forward neural networks (named also multi-layered perceptrons), therefore the following description is focused on such a class of networks. Since we deal with a scalar equation, a neural network will be a function $w:\mathbb{R}^n\rightarrow\mathbb{R}$ defined as follows: for any $\boldsymbol{x}\in\mathbb{R}^n$, the output $w(\boldsymbol{x})$ is computed via the chain of assignments
\begin{equation} \label{nn_main_formula}
  \begin{aligned}
  &\boldsymbol{x}_0=\boldsymbol{x}, \\
 &\boldsymbol{x}_\ell = \rho({\mathbf A}_\ell \boldsymbol{x}_{\ell-1} + {\mathbf b}_\ell), \hspace{1cm} \ell = 1,...,L-1, \\
 &w(\boldsymbol{x}) = {\mathbf A}_{L} \boldsymbol{x}_{L-1} + b_L\,.
  \end{aligned}
\end{equation}
Here, ${\mathbf A}_\ell\in\mathbb{R}^{N_\ell\times N_{\ell-1}}$ and ${\mathbf b}_\ell\in\mathbb{R}^{N_\ell}$, $\ell=1,...,L$, are matrices and vectors that store the network weights (with $N_0=n$ and $N_L=1$); furthermore, $L$ is the number of layers, whereas $\rho$ is the (nonlinear) activation function which acts component-wise (i.e. $\rho(\mathbf y)=\left[\rho(y_1),...,\rho(y_{n_y})\right]$ for any vector $\mathbf y\in\mathbb{R}^{n_y}$).  It can be noted from equation \eqref{nn_main_formula}, that if $\rho\in \Cr^k(\mathbb{R})$, then $w$ inherits the same regularity because it can be seen as a composition of functions belonging to $\Cr^k(\mathbb{R})$. Popular choices include the ReLU ($k=0$) and RePU ($k>0$ finite) functions, as well as the hyperbolic tangent ($k=\infty$) if one wants to exploit the maximum of regularity in the solution of interest.

The neural network structure ${\cal N\!N}$ is identified by fixing the number of layers $L$, the integers $N_\ell$ and the activation function $\rho$. The entire set of weights that parametrize the network can be logically organized into a single vector ${\mathbf w} \in \mathbb{R}^N$. Thus, the neural network structure ${\cal N\!N}$ induces a mapping
\begin{equation}\label{eq:generalNN}
{\cal F}^{\cal N\!N} : \mathbb{R}^N \to \Cr^\infty(\bar{\Omega})\,, \qquad {\mathbf w} \mapsto {\cal F}^{\cal N\!N} ({\mathbf w}) = w \,, \quad \text{where  } w = w(\boldsymbol{x},{\mathbf w})\,.
\end{equation}
It is convenient to define the manifold 
$$
U^{\cal N\!N} \subset U\,, \qquad U^{\cal N\!N} = {\cal F}^{\cal N\!N}(\mathbb{R}^N)
$$
containing all functions that can be generated by the neural network structure ${\cal N\!N}$.

\subsection{The VPINN discretization} \label{sec:sub_discretization}
We aim at approximating the solution of Problem \eqref{eq:model-pb} by a generalized Petrov-Galerkin strategy. To this end,  let us introduce  a conforming, shape-regular triangulation ${\cal T}_h= \{ E \}$ of $\bar{\Omega}$ with meshsize $h>0$ and, for a fixed integer $k_\text{test} \geq 1$, let $V_h \subset V$ be the linear subspace formed by the functions which are  piecewise polynomials of degree $k_\text{test}$ over the triangulation ${\cal T}_h$. Furthermore, let us introduce computable approximations of the forms $a$ and $F$ by numerical quadratures. Precisely, for any $E \in {\cal T}_h$, let $\{(\xi^E_\iota,\omega^E_\iota) : \iota \in I^E\}$ be the nodes and weights of a quadrature formula of precision 
\begin{equation}\label{eq:q-constraint}
q \geq 2k_\text{test}  
\end{equation}
on $E$. Assume that $\Gamma_N$ is the union of a collection $\partial {\cal T}_h (\Gamma_N)$ of edges of elements of ${\cal T}_h$; for any such edge $e$, let $\{(\xi^e_\iota,\omega^e_\iota) : \iota \in I^e\}$ be the nodes and weights of a quadrature formula of precision $q$ on $e$. Then, assuming that all the data $\mu$, $\boldsymbol{\beta}$, $\sigma$, $f$, $\psi$ are continuous in each element of the triangulation, we define the approximate forms
\begin{equation}\label{eq:def-ah}
a_h(w,v)= \sum_{E \in {\cal T}_h} \sum_{\iota \in I^E} [\mu \nabla w \cdot \nabla v + \boldsymbol{\beta}\cdot \nabla w \, v + \sigma w v](\xi^E_\iota) \,\omega^E_\iota\,, 
\end{equation}
\begin{equation}\label{eq:def-Fh}
F_h(v) =  \sum_{E \in {\cal T}_h}  \sum_{\iota \in I^E} [ f v](\xi^E_\iota) \,\omega^E_\iota \quad +  \sum_{e \in \partial{\cal T}_h(\Gamma_N)} \sum_{\iota \in I^e} [ \psi v](\xi^e_\iota) \, \omega^e_\iota \,.
\end{equation}

With these ingredients at hand, we would like to approximate the solution of Problem \eqref{eq:model-pb-var-2} by some $u^{\cal N\!N} \in U^{\cal N\!N}$ satisfying
\begin{equation}\label{eq:PGproblem}
a_h(Bu^{\cal N\!N},v_h)=F_h(v_h) \qquad \forall v_h \in V_h\,. 
\end{equation}
Such a problem might be ill-posed when,  for computational efficiency, the dimension of the test space $V_h$ is chosen smaller than the dimension of the manifold $U^{\cal N\!N}$. In this situation, we get an under-determined problem, with obvious difficulties in deriving stability estimates on some norms of the function $Bu^{\cal N\!N}$.
Actually, Problem \eqref{eq:PGproblem} with zero data (i.e., zero $f$, $g$, $\psi$) could admit non-zero solutions (see Section \ref{sect:inf-sup-section}).

To avoid these difficulties, we adopt the strategy of applying a projection (indeed, an interpolation) to the function $Bu^{\cal N\!N}$, mapping it into a finite dimensional space of dimension comparable to that of $V_h$, and we \cblue{limit} ourselves with estimating some norm of this projection.

To be precise, let us introduce a conforming, shape-regular partition ${\cal T}_H=\{G\}$ of $\bar{\Omega}$, which is equal to or coarser than ${\cal T}_h$ (i.e., each element $E \in {\cal T}_h$ is contained in an element $G \in {\cal T}_H$) but compatible with ${\cal T}_h$ (i.e., its meshsize $H>0$ satisfies $H\lesssim h$). Let the integer $\kt \geq 1$ be defined by the condition 
\begin{equation}\label{eq:def-ktrial}
\kt+k_\text{test} =q + 2 \,.
\end{equation}
Let $U_H \subset U$ be the linear subspace formed by the functions which are  piecewise polynomials of degree $\kt$ over the triangulation ${\cal T}_H$, and let $U_{H,0} = U_H \cap V$ be the subspace of $U_H$ formed by the functions vanishing on $\Gamma_D$. Finally, let ${\cal I}_H : \Cr^0(\bar{\Omega}) \to U_H$ be an interpolation operator, satisfying the condition
${\cal I}_H : \Cr^0(\bar{\Omega})\cap V \to U_{H,0}$ as well as the 
following approximation properties: for all $v \in \Hr^{k+1}(\Omega)$, $1 \leq k \leq \kt$,
\begin{equation}\label{eq:error-approx-v}
\vert v - {\cal I}_H v \vert _{\ell, G} \lesssim H^{k+1-\ell} \vert v\vert _{k+1,G}\,, \qquad 0 \leq \ell \leq k+1\,, \quad \forall G \in {\cal T}_H\,.
\end{equation}

In this framework, assuming the lifting $\bar{u}$ to be continuous in $\bar{\Omega}$, we replace the target equations  \eqref{eq:PGproblem} by the following ones:
\begin{equation}\label{eq:PGproblem-2}
a_h({\cal I}_H Bu^{\cal N\!N},v_h)=F_h(v_h) \qquad \forall v_h \in V_h\,. 
\end{equation}
In order to handle this problem with the neural network, let us introduce a basis in $V_h$, say $V_h = \text{span}\{\varphi_i : i\in I_h\}$, and for any $w$ smooth enough let us define the residuals
\begin{equation}\label{eq:residuals}
r_{h,i}(w)=F_h(\varphi_i)-a_h({\cal I}_H Bw,\varphi_i)\,, \qquad i \in I_h\,,
\end{equation}
as well as the loss function
\begin{equation}\label{eq:loss-function}
R_h^2(w) = \sum_{i \in I_h} r_{h,i}^2(w) \, \gamma_i^{-1}\,, 
\end{equation}
where $\gamma_i >0$ are suitable weights. Then, we search for a global minimum of the loss function in $U^{\cal N\!N}$, i.e., we consider the following discretization of Problem \eqref{eq:model-pb-var-2}: {\it Find $u^{\cal N\!N} \in U^{\cal N\!N}$ such that}
\begin{equation}\label{eq:min-prob}
u^{\cal N\!N} \in \displaystyle{\text{arg}\!\!\!\!\min_{w \in U^{\cal N\!N}}}\, R_h^2(w) \,.
\end{equation}
Note that the solution $u^{\cal N\!N}$ may not be unique; however, a suitable choice of the space $U_H$ may lead to the control of the error $u-{\cal I}_H Bu^{\cal N\!N}$ in the $H^1$-norm, as we will see in the sequel.

\begin{remark}[Discretization without interpolation]\label{rem:no-interp}{\rm
For the sake of comparison, we will also consider the optimization problem in which no interpolation is applied to the neural network functions. In other words, the target equations are those in \eqref{eq:PGproblem}, which induce the following definition of loss function
\begin{equation}\label{eq:loss-function-wI}
\hat{R}_h^2(w) = \sum_{i \in I_h} \hat{r}_{h,i}^2(w) \, \gamma_i^{-1}\,, \qquad \text{with } \quad  \hat{r}_{h,i}(w)=F_h(\varphi_i)-a_h(Bw,\varphi_i) \,,
\end{equation}
and the following minimization problem:  {\it Find $\hat{u}^{\cal N\!N} \in U^{\cal N\!N}$ such that}
\begin{equation}\label{eq:min-prob-wI}
\hat{u}^{\cal N\!N} \in \displaystyle{\text{arg}\!\!\!\!\min_{w \in U^{\cal N\!N}}}\, \hat{R}_h^2(w) \,.
\end{equation}
Note that in this problem the triangulation ${\cal T}_H$ and the space $U_H$ play no role. Although we will not provide a rigorous error analysis for such discretization, it will be interesting to numerically compare the behaviour of the approaches (i.e., with or without interpolation). This will be done in Section \ref{sect:num-res}. 
}
\end{remark}

\section{A priori error estimates} \label{sec:error_estimates}

Let $u^{\cal N\!N} \in U^{\cal N\!N}$ be any solution of the minimization problem \eqref{eq:min-prob}; let us set 
\begin{equation}\label{eq:def-uHNN}
u_H^{{\cal N\!N}} = {\cal I}_H Bu^{\cal N\!N} \in U_H\,.
\end{equation}
Recalling the definition \eqref{eq:split-u} of the affine mapping $B$, it holds
\begin{equation}\label{eq:decomp-uHNN}
u_H^{{\cal N\!N}}  = \bar{u}_H + u_H^{{\cal N\!N},0}\,, \qquad \text{with} \quad \bar{u}_H = {\cal I}_H \bar{u} \quad \text{and} \quad u_{H,0}^{\cal N\!N} = {\cal I}_H (\Phi u^{\cal N\!N}) \in U_{H,0}\,;
\end{equation}
note that $\bar{u}_H$ is a discrete lifting in $U_H$ of the Dirichlet data $g$.
%

\smallskip
We aim at estimating the error between $u$ and $u_H^{\cal N\!N}$.
To accomplish this task, we need several definitions, assumptions, and technical results.

\begin{definition}[norm-equivalence]\label{def:norm-equiv-Vh}
Let us denote by $0< c_h \leq C_h$ the constants in the norm equivalence
\begin{equation}\label{eq:norm-equiv-Vh}
c_h \Vert v_h \Vert_{1, \Omega} \leq \Vert \boldsymbol{v} \Vert_\gamma \leq C_h \Vert v_h \Vert_{1, \Omega} \qquad \forall v_h \in V_h \,, 
\end{equation}
where $\boldsymbol{v} = (v_i)_{i \in I_h}$ is such that $v_h = \sum_{i \in I_h} v_i \varphi_i$, and $\Vert \boldsymbol{v} \Vert_\gamma = \left( \sum_{i \in I_h} v_i^2\gamma_i \right)^{1/2}$. 
\end{definition}

Next, we introduce the consistency errors due to numerical quadratures 
\begin{equation}\label{eq:error-a} 
E_h^a(w_H,v_h) = a(w_H,v_h)-a_h(w_H,v_h) \qquad \forall w_H \in U_H, \ \forall v_h \in V_h \,,
\end{equation}
\begin{equation}\label{eq:errorF} 
E_h^F(v_h) = F(v_h)-F_h(v_h) \qquad \forall v_h \in V_h \,,
\end{equation}
and we provide a bound on these errors. To this end, let us assume that the quadrature rules used in the elements in ${\cal T}_h$ are obtained by affine transformations from a quadrature rule $\{(\hat{\xi}_\iota,\hat{\omega}_\iota) : \iota \in \hat{I}\}$ on a reference element $\hat{E} \subset \mathbb{R}^n$; similarly, let us assume that the quadrature rules used in the edges on $\partial{\cal T}_h(\Gamma_N)$ are obtained by affine transformations from a quadrature rule $\{(\check{\xi}_\iota,\check{\omega}_\iota) : \iota \in \check{I}\}$ on a reference element $\check{e} \subset \mathbb{R}^{n-1}$.

\begin{assumption}[Data smoothness]\label{ass:data} 
Let us assume the following smoothness of data:
\begin{equation}\label{eq:data-smoothness}
\mu, \ \sigma, \ f \in \Wr^{k,\infty}(\Omega)\,, \qquad \boldsymbol{\beta} \in (\Wr^{k,\infty}(\Omega))^n\,, \qquad \psi \in \Wr^{k,\infty}(\Gamma_N)\,,
\end{equation}
where $k$ is an integer satisfying
\begin{equation}\label{eq:def-k}
1 \leq k \leq \kt=q+2-k_\text{test}\,.
\end{equation}
\end{assumption}
Consequently, let us introduce the following notation
\begin{eqnarray}
{\cal N}_k(\mu, \boldsymbol{\beta}, \sigma) &=& \Vert \mu \Vert_{\Wr^{k,\infty}(\Omega)} + \Vert \boldsymbol{\beta} \Vert_{(\Wr^{k,\infty}(\Omega))^n} +  \Vert \sigma \Vert_{\Wr^{k,\infty}(\Omega)} \,, \\[5pt]
{\cal N}_k(f, \psi) &=& \Vert f \Vert_{\Wr^{k,\infty}(\Omega)} +  \Vert \psi \Vert_{\Wr^{k,\infty}(\Gamma_N)} \,,  \\
\Vert w_H \Vert_{k, {\cal T}_H} &=& \left( \sum_{G \in {\cal T}_H} \Vert w_{H \vert G} \Vert_{\Hr^k(G)} \right)^{1/2}  \quad \forall w_H \in U_H\,. 
\end{eqnarray}

\begin{property}[approximation of the forms $a$ and $F$]\label{prop:error-aF}
Under Assumption \ref{ass:data}, it holds
\begin{equation}\label{eq:error-a-est} 
\vert E_h^a(w_H,v_h)\vert   \ \lesssim \  h^k {\cal N}_k(\mu, \boldsymbol{\beta}, \sigma)  \Vert w_H \Vert_{k, {\cal T}_H}   \Vert v_h \Vert_{1,\Omega}  
                \qquad \forall w_H \in U_H, \ \forall v_h \in V_h \,,
\end{equation}
\begin{equation}\label{eq:error-F} 
\vert E_h^F(v_h)\vert   \ \lesssim \  h^k {\cal N}_k(f,\psi)   \Vert v_h \Vert_{1,\Omega}     \qquad \forall v_h \in V_h \,,
\end{equation}
\end{property}
\proof Both estimates are classical in the theory of finite elements (see, e.g., \cite{ciarlet2002finite}). As far as \eqref{eq:error-a-est} is concerned, the standard proof given for the case in which the polynomial degree is the same for both arguments, i.e., $k=k_\text{test} \geq 1$ and $q=2(k-1)$, can be easily adapted to the present situation $k+k_\text{test} \leq q+2$. In this way, one gets $\vert E_h^a(w_H,v_h)\vert   \ \lesssim \  h^k {\cal N}_k(\mu, \boldsymbol{\beta}, \sigma)  \Vert w_H \Vert_{k, {\cal T}_h}   \Vert v_h \Vert_{1,\Omega}$, and one concludes by observing that
$\Vert w_H \Vert_{k, {\cal T}_h} = \Vert w_H \Vert_{k, {\cal T}_H}$ since  ${\cal T}_h$ is a refinement of  ${\cal T}_H$. 
\endproof

\medskip
Finally, we pose a fundamental assumption.

\begin{assumption}[inf-sup condition between $U_{H,0}$ and $V_h$]\label{ass:inf-sup}
The bilinear form $a$ satisfies an inf-sup condition with respect to the spaces $U_{H,0}$ and $V_h$, namely there exists a constant $\alpha_\star >0$, independent of the meshsizes $h$ and $H$, such that
\begin{equation}\label{eq:infsup}
\alpha_\star \Vert w_H \Vert_{1,\Omega} \leq \sup_{v_h \in V_h} \frac{a(w_H,v_h)}{\Vert v_h \Vert_{1,\Omega}} \qquad \forall w_H \in U_{H,0}\,.
\end{equation}
\end{assumption}


This assumption together with Property \ref{prop:error-aF} yields the following result.

\begin{proposition}[discrete inf-sup condition between $U_{H,0}$ and $V_h$]\label{cor:discrete-inf-sup}
 Under Assumptions \ref{ass:data} and \ref{ass:inf-sup}, for all $h \leq h_0$ small enough the bilinear form $a_h$ satisfies an inf-sup condition with respect to the spaces $U_{H,0}$ and $V_h$, namely there exists a constant $\tilde{\alpha}_\star >0$ such that
\begin{equation}\label{eq:infsup-h}
\tilde{\alpha}_\star \Vert w_H \Vert_{1,\Omega} \leq \sup_{v_h \in V_h} \frac{a_h(w_H,v_h)}{\Vert v_h \Vert_{1,\Omega}} \qquad \forall w_H \in U_{H,0}\,.
\end{equation}
\end{proposition}
\proof We have $a_h(w_H,v_h) = a(w_H,v_h) - E_h^a(w_H,v_h)$. Using the bound \eqref{eq:error-a-est} with $k=1$ and observing that $\Vert w_H \Vert_{1,{\cal T}_H}=\Vert w_H \Vert_{1,\Omega}$, one can find $h_0>0$ small enough such that, for all $h \leq h_0$, $\vert E_h^a(w_H,v_h)\vert    \leq \frac12 \alpha_\star \Vert w_H \Vert_{1,\Omega}   \Vert v_h \Vert_{1,\Omega}$, whence the result with  $\tilde{\alpha}_\star = \frac12 \alpha_\star$ \endproof

\bigskip
We are ready to estimate the error $\Vert u - u_H^{\cal N\!N} \Vert_{1,\Omega}$. Recalling the decomposition \eqref{eq:decomp-uHNN},  we use the triangle inequality
\begin{equation}\label{eq:bound1}
\Vert u - u_H^{\cal N\!N} \Vert_{1,\Omega} \leq \Vert u - u_H \Vert_{1,\Omega}  + \Vert u_H - u_H^{\cal N\!N} \Vert_{1,\Omega} \,,
\end{equation}
where $u_H$ is a suitable element in the affine subspace $\bar{u}_H + U_{H,0} \subset U_H$. Writing $u_H = \bar{u}_H + u_{H,0}$ with $u_{H,0} \in U_{H,0}$,
one has $u_H - u_H^{\cal N\!N} = u_{H,0}  - u_{H,0}^{\cal N\!N} \in U_{H,0}$; hence, we can apply \eqref{eq:infsup-h} to get
\begin{equation}\label{eq:bound2}
\Vert u_H - u_H^{\cal N\!N} \Vert_{1,\Omega} \leq \frac1{\tilde{\alpha}_\star} \sup_{v_h \in V_h} \frac{a_h(u_H - u_H^{\cal N\!N},v_h)}{\Vert v_h \Vert_{1,\Omega}} \,.
\end{equation}
Recalling the definitions \eqref{eq:error-a} and \eqref{eq:errorF}, it holds
\begin{eqnarray*}
a_h(u_H,v_h) &=& a(u_H,v_h) - E_h^a(u_H,v_h) \\
&=& a(u,v_h) - a(u-u_H,v_h) - E_h^a(u_H,v_h) \\
&=& F(v_h) - a(u-u_H,v_h) - E_h^a(u_H,v_h) \\
&=& F_h(v_h) +E_h^F(v_h)- a(u-u_H,v_h) - E_h^a(u_H,v_h) \,.
\end{eqnarray*}
Thus, the numerator in \eqref{eq:bound2} is given by
$$
a_h(u_H - u_H^{\cal N\!N},v_h) = F_h(v_h) - a_h(u_H^{\cal N\!N},v_h) - a(u-u_H,v_h) - E_h^a(u_H,v_h) +E_h^F(v_h)\,.
$$
On the other hand, recalling \eqref{eq:residuals} we have
\begin{equation}\label{eq:residual-repr}
F_h(v_h) - a_h(u_H^{\cal N\!N},v_h)=F_h(v_h) - a_h({\cal I}_H Bu^{\cal N\!N},v_h) = \sum_{i \in I_h} r_{h,i}(u^{\cal N\!N}) \, v_i \,,
\end{equation}
hence, by \eqref{eq:loss-function} and \eqref{eq:norm-equiv-Vh},
$$
\vert  F_h(v_h) - a_h(u_H^{\cal N\!N},v_h) \vert  \leq R_h(u^{\cal N\!N})  \Vert \boldsymbol{v} \Vert_\gamma \leq C_h R_h(u^{\cal N\!N}) \Vert v_h \Vert_{1, \Omega} \,.
$$
Using the bounds \eqref{eq:error-a-est}  and \eqref{eq:error-F}, we obtain the following inequality
\begin{equation}\label{eq:pre-lemma}
\begin{split}
\Vert u - u_H^{\cal N\!N} \Vert_{1,\Omega} & \lesssim \left( 1+\frac1{\tilde{\alpha}_\star} \right) \Bigl(\ \inf_{u_H \in \bar{u}_H +U_{H,0}}\left(  \Vert u - u_H \Vert_{1,\Omega} 
+ h^k {\cal N}_k(\mu, \boldsymbol{\beta}, \sigma)  \Vert u_H \Vert_{k, {\cal T}_H}  \right)  \\
&\qquad \qquad \qquad \qquad + C_h R_h(u^{\cal N\!N}) + h^k {\cal N}_k(f,\psi)  \  \Bigr)\,.
\end{split}
\end{equation}

From now on, we assume that $u \in \Hr^{k+1}(\Omega)$. Then, assumption \eqref{eq:error-approx-v} yields the inequalities
\begin{equation}\label{eq:error-approx-u-norm}
 \Vert u - {\cal I}_H u \Vert_{1,\Omega}  \lesssim H^k \vert u\vert _{k+1,\Omega} \lesssim h^k \vert u\vert _{k+1,\Omega}
\end{equation}
and
\begin{equation}\label{eq:bound-PiHu-k}
\Vert {\cal I}_H u \Vert_{k, {\cal T}_H} \leq \Vert u \Vert_{k,\Omega} + \Vert u-{\cal I}_H u \Vert_{k, {\cal T}_H} 
 \lesssim \Vert u \Vert_{k,\Omega} + H \vert u\vert _{k+1,\Omega} \lesssim \Vert u \Vert_{k+1,\Omega} \,.
\end{equation}
Choosing $u_H = {\cal I}_H u \in \bar{u}_H +U_{H,0}$ in \eqref{eq:pre-lemma} and using these estimates, we arrive at the following intermediate result, \cblue{which can be viewed as a mixed a priori/a posteriori error estimate.}
\begin{lemma}\label{lem:intermedio1}
Under the previous assumptions, it holds
$$
\Vert u - u_H^{\cal N\!N} \Vert_{1,\Omega}  \lesssim  h^k (\vert u\vert _{k+1,G} + {\cal N}_k(\mu, \boldsymbol{\beta}, \sigma) \Vert u \Vert_{k+1,\Omega} + {\cal N}_k(f,\psi)  ) + C_h R_h(u^{\cal N\!N}) \,.
$$
\end{lemma}

\medskip
Our next task will be bounding the term $R_h(u^{\cal N\!N})$. To this end, we use the minimality condition \eqref{eq:min-prob} to get
\begin{equation}\label{eq:u-minimizer}
R_h(u^{\cal N\!N}) \leq R_h(w^{\cal N\!N}) \qquad \forall w^{\cal N\!N} \in U^{\cal N\!N} \,.
\end{equation}
On the other hand, since $R_h(w^{\cal N\!N})$ is a weighted $\ell_2$-norm in $\mathbb{R}^{\vert I_h\vert }$, we can write
$$
R_h(w^{\cal N\!N}) = \sup_{{\mathbf z} \in \mathbb{R}^{\vert I_h\vert }} \frac1{\Vert {\mathbf z} \Vert_\gamma} \sum_{i \in I_h} r_{h,i}(w^{\cal N\!N}) z_i \,,
$$
where, similarly to \eqref{eq:residual-repr},
$$
\sum_{i \in I_h} r_{h,i}(w^{\cal N\!N}) z_i  = F_h(z_h) - a_h({\cal I}_H Bw^{\cal N\!N}, z_h) \qquad \text{with \ } z_h = \sum_{i \in I_h}  z_i \varphi_i \in V_h \,.
$$
For convenience, in analogy with \eqref{eq:def-uHNN}, let us set
\begin{equation}\label{eq:def-wHNN}
w_H^{{\cal N\!N}} = {\cal I}_H Bw^{\cal N\!N} \in U_H\,.
\end{equation}
Thus, recalling \eqref{eq:norm-equiv-Vh}, we obtain
\begin{equation}\label{eq:bound3}
R_h(w^{\cal N\!N}) \leq  \frac1{c_h} \sup_{z_h \in V_h} \frac{F_h(z_h) - a_h(w_H^{\cal N\!N}, z_h)}{\Vert z_h \Vert_{1, \Omega}} \qquad \forall  w^{\cal N\!N} \in U^{\cal N\!N} \,.
\end{equation}
The numerator can be manipulated as above, using
$$
F_h(z_h) = F(z_h) - E_h^F(z_h) =  a(u,z_h)  - E_h^F(z_h)  
$$
and
$$
a_h(w_H^{\cal N\!N}, z_h) = a(w_H^{\cal N\!N}, z_h) - E_h^a(w_H^{\cal N\!N},z_h) \,,
$$
whence, using once more Property \ref{prop:error-aF}, we get
\begin{equation}\label{eq:boundR_h-intermed}
R_h(w^{\cal N\!N}) \lesssim  \frac1{c_h} \left( \Vert u - w_H^{\cal N\!N} \Vert_{1,\Omega} 
                 + h^k {\cal N}_k(\mu, \boldsymbol{\beta}, \sigma)  \Vert w_H^{\cal N\!N}  \Vert_{k, {\cal T}_H} + h^k {\cal N}_k(f,\psi)  \right)\,.
\end{equation}

In order to bound the terms containing $w_H^{\cal N\!N}$, we introduce the quantity
\begin{equation}\label{eq:def-eNN}
e^{{\cal N\!N}} = u-Bw^{\cal N\!N}\,,
\end{equation}
which, recalling the definitions \eqref{eq:split-u} and \eqref{eq:def-B}, can be written as 
\begin{equation}\label{eq:def-eNN-1}
e^{{\cal N\!N}} = \Phi(\tilde{u}-w^{\cal N\!N}) \,,
\end{equation}
and we formulate a final assumption.

\begin{assumption}[smoothness of the solution and the neural network manifold]\label{ass:sol+NN} 
The solution $u$ can be represented as in \eqref{eq:split-u} with
\begin{equation}\label{eq:ass-reg-tuPhi}
\tilde{u} \in \Hr^{k+1}(\Omega) \qquad \text{and} \qquad \Phi \in \Wr^{k+1,\infty}(\Omega)
\end{equation}
for $k$ satisfying \eqref{eq:def-k}.  Furthermore, the manifold formed by the neural network functions satisfies the smoothness condition
\begin{equation}\label{eq:ass-reg-UNN}
U^{\cal N\!N} \subset \Hr^2(\Omega)\,.
\end{equation}
\end{assumption}
Note that \eqref{eq:ass-reg-tuPhi} implies in particular $u \in \Hr^{k+1}(\Omega)$ with $\Vert u \Vert_{k+1,\Omega} \lesssim \Vert \tilde{u} \Vert_{k+1,\Omega}\,\Vert \Phi \Vert_{k+1,\infty,\Omega}$; on the other hand,
\eqref{eq:ass-reg-UNN} implies $e^{\cal N\!N} \in \Hr^2(\Omega)$. (We refer to Remark \ref{rem:relu} for another set of assumptions on the neural network.)

Recalling \eqref{eq:def-wHNN} and using the identity
\begin{equation}\label{eq:decomp-error}
u-w_H^{\cal N\!N}  = (u-{\cal I}_H u) + {\cal I}_H e^{\cal N\!N} = (u-{\cal I}_H u) - e^{\cal N\!N} + (I-{\cal I}_H) e^{\cal N\!N} \,,
\end{equation}
we can write 
\begin{equation}\label{eq:bound-u-wHNN}
\Vert u-w_H^{\cal N\!N} \Vert_{1,\Omega}  \lesssim  \Vert u-{\cal I}_H u \Vert_{1,\Omega}  + \Vert e^{\cal N\!N} \Vert_{1,\Omega} + H \vert e^{\cal N\!N}\vert _{2,\Omega} 
\end{equation}
and, using a standard inverse inequality in $\mathbb{P}_k(G)$ for any $G \in {\cal T}_H$,
\begin{equation}\label{eq:bound-wHNN-k}
\begin{split}
\Vert w_H^{\cal N\!N}  \Vert_{k, {\cal T}_H}  &\lesssim  \Vert {\cal I}_H u \Vert_{k, {\cal T}_H} + H^{1-k} \Vert {\cal I}_H e^{\cal N\!N} \Vert_{1,\Omega} \\
&\lesssim \Vert {\cal I}_H u \Vert_{k, {\cal T}_H} + H^{1-k} \left( \Vert  e^{\cal N\!N} \Vert_{1,\Omega} + H  \vert e^{\cal N\!N}\vert _{2,\Omega} \right) \,.
\end{split}
\end{equation}
Keeping into account \eqref{eq:error-approx-u-norm} and \eqref{eq:bound-PiHu-k},  in order to conclude we need to identify a function $\tilde{w}^{\cal N\!N} \in U^{\cal N\!N}$ for which a bound of the type  
\begin{equation}\label{eq:bound-eNN}
\vert  e^{\cal N\!N}\vert _{m,\Omega} \lesssim \vert  \tilde{u}-\tilde{w}^{\cal N\!N}\vert _{m,\Omega} \lesssim  H^{k+1-m} \vert \tilde{u}\vert _{k+1,\Omega}
\end{equation}
holds true for $m=1,2$. The existence of such a function is guaranteed by one of the available results on the approximation of functions in Sobolev spaces by neural networks (see \cite[Theorem 5.1, Remark 5.2]{deryck2021approximation}; see also \cite{opschoor2021exponential}), provided the number of layers $L$ and the widths of the layers in the chosen ${\cal N\!N}$ satisfy suitable conditions depending on the target accuracy (hence, in our case depending on $H^k$). 
\cblue{Indeed, suppose one is interested in using meshes with meshsize as small as $H_{\min }$ in the domain $\Omega$ (here assumed to satisfy $\Omega\subset [0,1]^n$ for the sake of simplicity), 
and let $N\in\mathbb N$ be such that
\begin{equation}\label{eq:n_error_for_tanh}
3^n(1+\delta)(2(m+1))^{3m} \max\{R^m,\text{ln}^m(\beta N^{k+n+3})\}\dfrac{C(n,m,k,\tilde u)}{N^{k+1-m}} \le  H_{\min }^{k+1-m} \vert \tilde{u}\vert _{k+1,\Omega} \,,
\end{equation}
where $\delta$, $R$, $\beta$ and $C$ are constants not depending on $N$ defined in \cite{deryck2021approximation}. Then, a function $\tilde{w}^{\cal N\!N}$ exists which fulfils \eqref{eq:bound-eNN} and is represented as a neural network with the hyperbolic tangent as activation function and two hidden layers with $N_1$ and $N_2$ neurons respectively, satisfying
\begin{equation}\label{eq:n1_n2_bounds}
N_1 \le 3 \left\lceil \frac {k+1}2\right\rceil \left\vert P_{k,n+1}\right\vert + n(N-1),
\hspace{0.3cm}
N_2 \le 3 \left\lceil \frac{n+2}2 \right\rceil \left\vert P_{n+1,n+1}\right\vert N^n,
\end{equation}
where 
\[
\left\vert P_{a,b}\right\vert = \left( a+b-1 \atop a \right),\hspace{1cm} \forall a,b\in\mathbb N, \ b\ge2.
\] 
}

Substituting \eqref{eq:bound-eNN} into \eqref{eq:bound-u-wHNN} and \eqref{eq:bound-wHNN-k}, and using inequalities \eqref{eq:u-minimizer} and  \eqref{eq:boundR_h-intermed}, we arrive at the following bound on the loss $R_h(u^{\cal N\!N})$. 
\begin{lemma}\label{lem:intermedio2}
Under the previous assumptions, it holds
$$
R_h(u^{\cal N\!N}) \lesssim  \frac1{c_h} \left( H^k  \vert u\vert _{k+1,\Omega}  + H^k  \vert \tilde{u}\vert _{k+1,\Omega} + h^k {\cal N}_k(\mu, \boldsymbol{\beta}, \sigma)  \Vert \tilde{u} \Vert _{k+1,\Omega}  + h^k {\cal N}_k(f,\psi)  \right)\,.
$$
\end{lemma}
\cblue{We remark that such a bound, when the involved neural network is comprised of at least two hidden layers and is such that there exists $N$ satisfying both \eqref{eq:n_error_for_tanh} and \eqref{eq:n1_n2_bounds}, does not depend on the network hyperparameters.}

Concatenating Lemmas \ref{lem:intermedio1} and \ref{lem:intermedio2}, and using once more $H \lesssim h$, we obtain the following a priori error estimate for the solution of Problem \eqref{eq:min-prob}.
\begin{theorem}[a priori error estimate]\label{teo:a-priori bound}
Let $u_H^{\cal N\!N} \in U_H$ be defined by \eqref{eq:def-uHNN}. Under Assumptions \ref{ass:data}, \ref{ass:inf-sup} and \ref{ass:sol+NN}, for $h$ sufficiently small it holds
\begin{equation}\label{eq:a-priori bound}
\begin{split}
\Vert u - u_H^{\cal N\!N} \Vert_{1,\Omega} &\lesssim \left( 1+\frac{C_h}{c_h} \right)  h^k \big[ (1+{\cal N}_k(\mu, \boldsymbol{\beta}, \sigma) ) \Vert \tilde{u} \Vert _{k+1,\Omega} + {\cal N}_k(f,\psi) \big] \,.
\end{split}
\end{equation}
\end{theorem}

\begin{remark}[on the equivalence constants $c_h, C_h$]\label{rem:equiv-ch}{\rm
If a classical Lagrange basis is used in \eqref{eq:residuals}, and the triangulation ${\cal T}_h$ is quasi-uniform, then for constants weights $\gamma_i=1$ one has $c_h \simeq h^{1-d/2}$ and
$C_h \simeq h^{-d/2}$, whence ${\frac{C_h}{c_h}} \simeq h^{-1}$.  On the other hand, if a hierarchical basis is used instead, then $c_h \simeq C_h \simeq 1$, hence, ${\frac{C_h}{c_h}} \simeq 1$ in dimension $d=1$, whereas $c_h \simeq \vert \log h\vert ^{-1}$, $C_h \simeq 1$, hence, $\frac{C_h}{c_h} \simeq \vert \log h\vert $ in dimension $d=2$. 

Thus, the presence of the ratio $\frac{C_h}{c_h}$ in \eqref{eq:a-priori bound}, which originates from the control of the loss function, makes this estimate sub-optimal. \cblue{However, our numerical experiments in Sect. \ref{sect:error_vs_h}  indicate that this adverse effect is not seen in practice. The reason may be related to the decay of the loss function $R_h(u^{\cal N\!N})$, which is significantly faster than the decay of the approximation error when $h$ is reduced, thereby compensating for the growth of ratio. See Remark \ref{rem:equiv-ch-bis}. }\qquad $\square$ }
\end{remark}

%

\begin{remark}[low-regularity ${\cal N\!N}$]\label{rem:relu}{\rm
When the condition $U^{\cal N\!N} \subset \Hr^2(\Omega)$ fails to be satisfied, as for the ReLU activation function, we may provide a different set of assumptions which still lead to an $O(h^k)$-error estimate as in Theorem \ref{teo:a-priori bound}. Precisely, we may assume that $\tilde{u} \in \Wr^{k+1,\infty}(\Omega)$ and $U^{\cal N\!N} \subset \Wr^{1,\infty}(\Omega)$. Then, referring to the first equality in \eqref{eq:decomp-error}, one has
$$
\Vert u-w_H^{\cal N\!N} \Vert_{1,\Omega} \leq  \Vert (u-{\cal I}_H u) \Vert_{1,\Omega} + \Vert {\cal I}_H e^{\cal N\!N} \Vert_{1,\Omega} \lesssim
\Vert (u-{\cal I}_H u) \Vert_{1,\Omega} + H^{-1}\Vert {\cal I}_H e^{\cal N\!N} \Vert_{0,\Omega} \,, 
$$
with
\begin{equation*}
\begin{split}
\Vert {\cal I}_H e^{\cal N\!N} \Vert_{0,\Omega}^2 &= \sum_{G \in {\cal T}_H} \Vert {\cal I}_H e^{\cal N\!N} \Vert_{0,G}^2 
\leq \sum_{G \in {\cal T}_H} \Vert {\cal I}_H e^{\cal N\!N} \Vert_{\Lr^\infty(G)}^2 \vert G\vert  \\
&\lesssim \sum_{G \in {\cal T}_H} \Vert e^{\cal N\!N} \Vert_{\Lr^\infty(G)}^2 \vert G\vert \lesssim H^{d} \Vert e^{\cal N\!N} \Vert_{\Lr^\infty(\Omega)}^2 \,.
\end{split}
\end{equation*}
The conclusion easily follows if $\tilde{w}^{\cal N\!N}$ is chosen to satisfy the error bound $\Vert  \tilde{u}-\tilde{w}^{\cal N\!N}\Vert _{\Lr^\infty(\Omega)} \lesssim  H^{k+1} \vert \tilde{u}\vert _{\Wr^{k+1,\infty}(\Omega)}$\,, which is possible according to the results in \cite{guhring2020error}, \cite{opschoor2020deep}. \quad $\square$
\endproof }
\end{remark}

%% file: implementation.tex
\section{Implementation issues}\label{sec:implementation}

As specified in Section \ref{sec:networks}, we use a fully-connected feed-forward neural network architecture, which is fixed and depends neither on the PDE nor on its discretization. 
For each simulation, we initialize the neural network with a completely new set of weights, which is important to show that our results are not initialization dependent. The activation function is the hyperbolic tangent. It has been proven in \cite{deryck2021approximation} that such neural networks with two hidden layers enjoy exponential converge properties with respect to the number of weights. Nevertheless, in order to simplify the training and enrich the space in which we seek the numerical solution, we consider five layers (namely, $L=5$ with the notation of Section \ref{sec:networks}) with 50 neurons each ($N_\ell=50,$ $\ell=1,...,5$). We also highlight that it is always possible to approximate the identity function with a neural network with a single layer with just one neuron. The best approximation obtainable with a neural network with more than two layers is thus more accurate than the one computable with a neural network with just two layers, because, in the worst possible case, the $L$-layers neural network can be obtained by combining an $(L-2)$-layers identity neural network with a suitable 2-layers neural network. \cblue{Numerical tests have been performed to investigate the influence of the activation function on the model accuracy; we observed that all the commonly used activation functions led to equivalent results, thus we omit such a comparison from the present work. We remark that the hyperparameters $L=5$ and $N_\ell=50$ have been chosen to obtain a neural network \cblue{sufficiently large} to satisfy condition \eqref{eq:bound-eNN} on all grids used in our experiments. Numerical evidence that the neural network best approximation error is negligible when compared with other sources of error is presented in Sect. \ref{sect:error_vs_dimension}.}

In order to compute the results shown in Section \ref{sect:num-res}, we compared various state of the art optimizers to find the most efficient way to minimize the loss function. We observe that most of the momentum-based first-order methods have similar performances (the presented results are computed with the \cblue{ADAM optimizer \cite{kingma2014adam}}), but it is convenient to use a learning rate scheduler in order to reduce the learning rate during the process. We tested both cyclical learning rate schedulers and exponential learning rate schedulers \cite{smith2017cyclical}, the differences were very subtle and we thus choose to adopt the most common exponential learning rate scheduler and decided not to report images about such a comparison. To further reduce the loss function we also use the BFGS method and its limited-memory version: the L-BFGS method \cite{wright1999numerical}. 

The Dirichlet boundary conditions are imposed via the mapping $B$ defined in \eqref{eq:def-B}. The construction of the function $\Phi$ is particularly simple when $\Omega$ is a convex polygon, since in this case $\Phi$ can be defined as the product of the linear polynomials which vanish on each Dirichlet edge; this is precisely how we define $\Phi$ in the numerical examples discussed in the next session. \cblue{In other geometries, one can build $\Phi$ either as described in \cite{sukumar2021exact}, or by using a level-set method, or even by training an auxiliary neural network to (approximately) vanish on $\Gamma_D$. Similarly, in order to obtain an analytical expression of the extension $\overline u$ of the Dirichlet data $g$, one can train another neural network to (approximately) match the values of $g$ on $\Gamma_D$ or use a data transfinite interpolation \cite{RVACHEV2001195}.}

\subsection{VPINN efficiency and the inf-sup condition}\label{sec:inf-sup-advantages}
In Sect. \ref{sec:discretization} we introduced a discretization method in which the loss function is built by a piecewise polynomial interpolation of the neural network; on the other hand, we also mentioned in Remark \ref{rem:no-interp} the possibility of building the loss function directly from the (non-interpolated) neural network. 

From the theoretical point of view, only the former approach can be considered mathematically reliable, since the error control is based on the validity of an inf-sup condition, as detailed in Sect. \ref{sec:error_estimates}. On the contrary, if the neural network is used without interpolation, one usually gets an under-determined system, for which the error control may be problematic. In fact, for instance, the discrete solution with zero data may not be identically zero, as documented in Sect. \ref{sect:inf-sup-section}, which rules out uniqueness.
Nonetheless, there is empirical evidence  (see, e.g., \cite{zhang2020physics, sahli2020physics, ji2020stiff, wight2020solving}) that non-interpolated neural networks do succeed in computing accurate solutions even in complex scenarios. Actually, in the next section we will provide numerical evidence that the two approaches are always (in the considered cases) comparable in terms of rate of convergence and, when the solution is regular, smaller errors are obtained minimizing the same loss function without interpolation.

From the computational point of view, the two approaches have comparable advantages and disadvantages. Let us first consider non-interpolated VPINNs. The corresponding loss functions can be more easily implemented thanks to the existing deep-learning frameworks, which allow the direct computation of neural network derivatives via automatic differentiation \cite{baydin2018automatic}. One only needs to generate a mesh and the corresponding test functions, associate a quadrature rule with each element and assemble all the tensors required to efficiently compute the loss function. The main difference with the interpolated neural network approach is that, in the latter, the interpolation matrices have to be assembled too \cblue{(see Appendix  \ref{sec:Ih_construction} for a detailed description of the construction of the interpolation operators)}, while automatic differentiation is not required. Depending on the problem at hand, this may be an advantage or not. Indeed, the interpolation matrices assembly may be tricky but, using them, all derivatives can be efficiently computed by matrix-vector multiplications that are much cheaper than the entire automatic differentiation procedure, especially when higher order derivatives are required. Therefore, for fast-converging optimization processes, a non-interpolated neural network approach may be efficient and can be more easily implemented; otherwise, each optimization step may be much more expensive than the analogous operation performed with an interpolated neural network. Furthermore, we observed that the training phase is faster when the neural network is interpolated because the procedure converges in fewer steps. This is probably related to the fact that the solution is sought in a significantly smaller space, that can be more easily explored during the training phase.

%% file: numerical_results.tex
\section{Numerical results}\label{sect:num-res}
\begin{figure}[t!]
  \centering\hspace{-2cm}
  \begin{tabular}[c]{cc}
    \multirow{2}{*}[30pt]{
    \begin{subfigure}{0.60\textwidth}
     \hspace{-1.cm} \includegraphics[width=\textwidth]{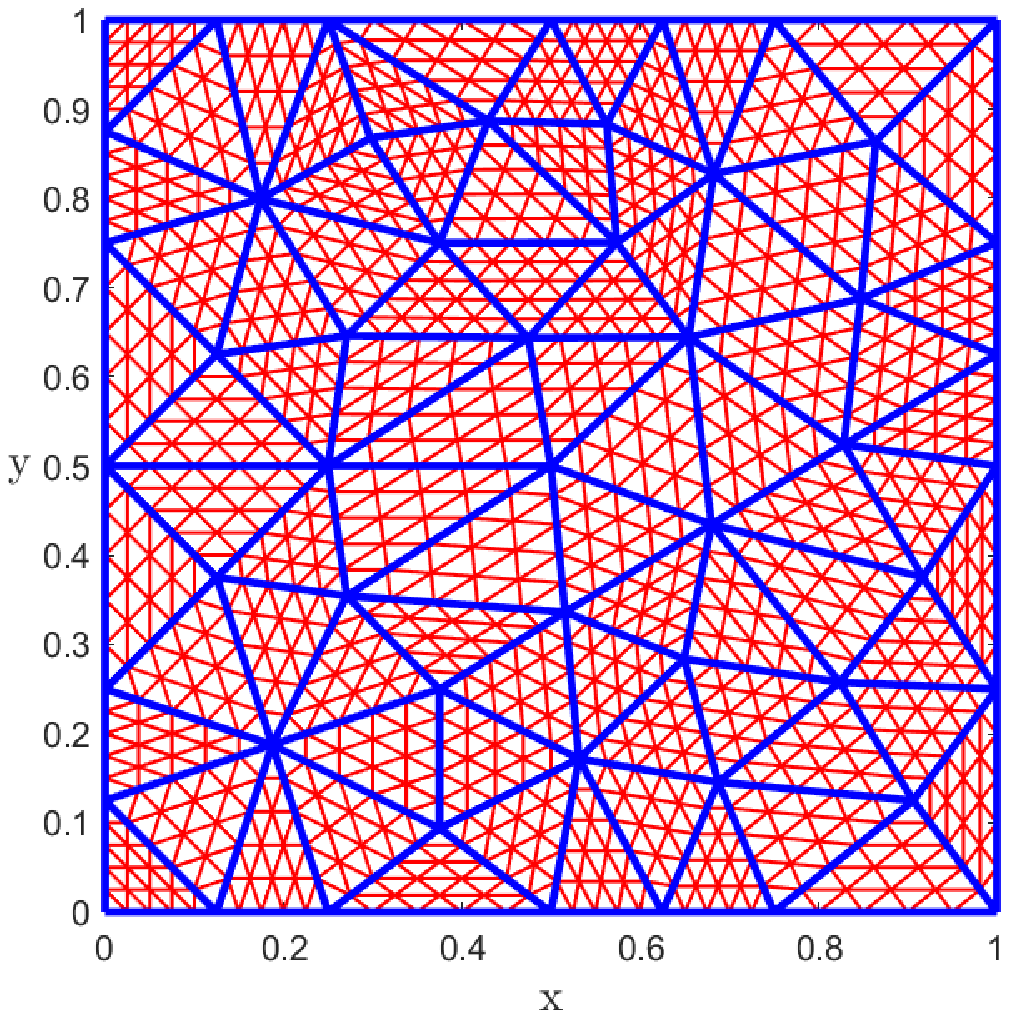} 
      \caption{ Delaunay meshes\\ $\mathcal T_H$ (blue lines) and $\mathcal T_h$ (red lines).\\ Refinement corresponding to $k_{\text{int}}=5$.}
      \label{fig:mesh_TH_Th}
    \end{subfigure}
}&\vspace{0.2cm}
   \begin{subfigure}[c]{0.13\textwidth}
      \includegraphics[width=\textwidth]{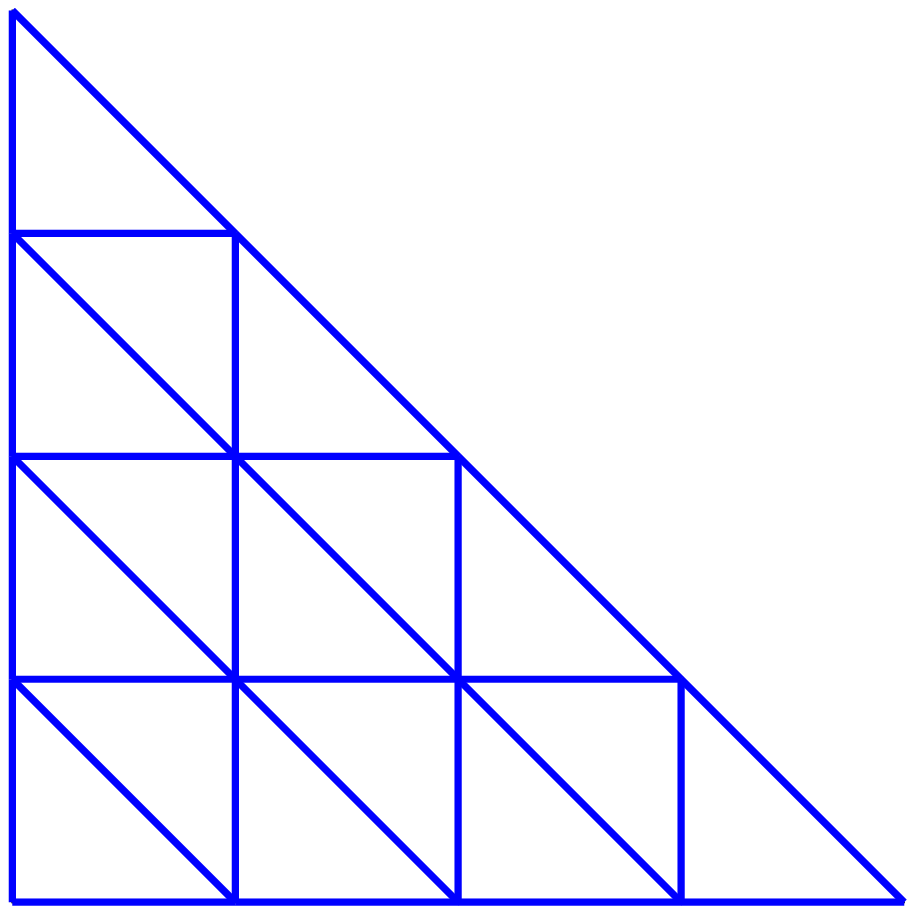} 
      \caption{\hspace{-1.6cm}\mbox{Elemental refinement for $k_{\text{int}}=4$.}}
      \label{fig:ktrial4}
    \end{subfigure}\\
    &\vspace{0.2cm}
        \begin{subfigure}[c]{0.13\textwidth}
      \includegraphics[width=\textwidth]{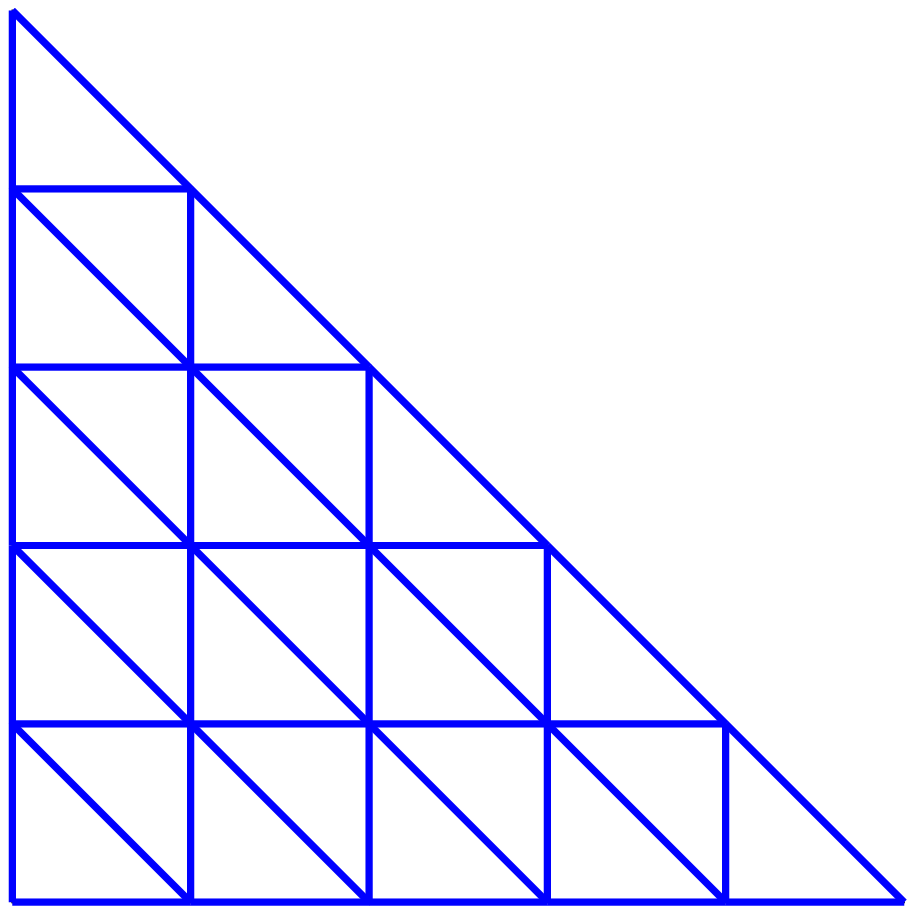} 
      \caption{\hspace{-1.6cm}\mbox{Elemental refinement for $k_{\text{int}}=5$.}}
      \label{fig:ktrial5}
    \end{subfigure}\\
    &\vspace{0.2cm}
   \begin{subfigure}[c]{0.13\textwidth}
      \includegraphics[ width=\textwidth]{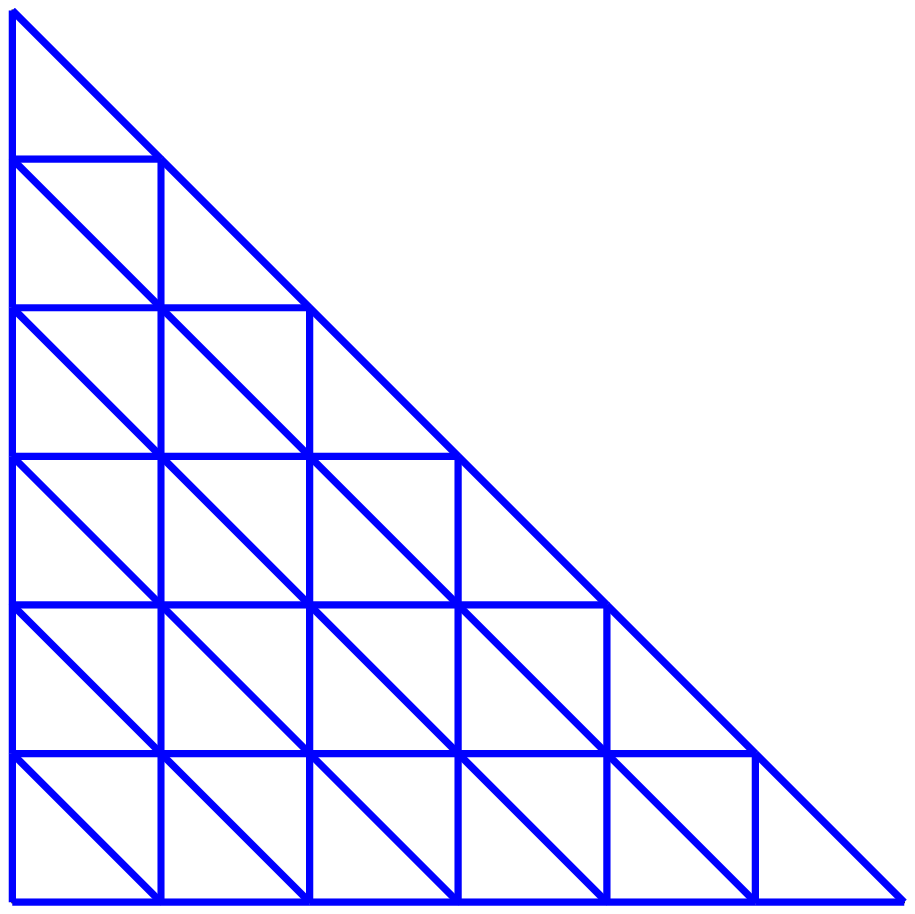} 
     \caption{ \hspace{-1.6cm}\mbox{Elemental refinement for $k_{\text{int}}=6$.}}
      \label{fig:ktrial6}
    \end{subfigure}\\
  \end{tabular}    
  \caption{One of the meshes used in Sect. \ref{sect:error_vs_h} and elemental refinements chosen to obtain $\mathcal T_h$ from $\mathcal T_H$.}
  \label{fig:mesh_refinement}
\end{figure}

\cblue{In this section we present several numerical results concerning the VPINN discretization of Problem \eqref{eq:model-pb} in the square $\Omega=(0,1)^2$. We will vary the coefficients of the operator, the boundary conditions and the smoothness of the exact solution. For each test case, we vary the degree $k_\text{test}$ of the test functions, the order $q$ of the quadrature rule and, correspondingly, we choose the polynomial degree $k$ of the interpolating functions as $k=\kt=q+2-k_\text{test}$, according to \eqref{eq:def-k}. We only report results obtained with Gaussian rules, as Newton-Cotes formulas of the same order give comparable errors \cblue{(see \cite{BeCaPi2021} for a larger set of numerical experiments about this and other comparisons)}.

The theoretical results in Sect. \ref{sec:error_estimates} suggest that it is convenient to maintain $k_\text{test}$ as low as possible; consequently, we only use piecewise linear ($k_\text{test}=1$) or piecewise quadratic ($k_\text{test}=2$) test functions. Recalling condition \eqref{eq:q-constraint}, we thus choose $q=3$ or $q=5$ if $k_\text{test}=1$, and $q=5$ if $k_\text{test}=2$.

The triangulations we use are generic Delaunay triangular meshes. In order to satisfy the discrete inf-sup condition, we choose $\mathcal T_H$ and $\mathcal T_h$ as nested meshes whose meshsizes satisfy $H=\kt h$. A pair $(\mathcal T_H, \mathcal T_h)$ of used meshes is represented in Fig. \ref{fig:mesh_refinement}, together with the elemental refinement corresponding to $k_{\text{int}}=4$, $k_{\text{int}}=5$ and $k_{\text{int}}=6$.
}

\cblue{\subsection{Error decays} \label{sect:error_vs_h}
Hereafter,  we empirically confirm, with numerical experiments, the a priori error estimate established in Sect.  \ref{sec:error_estimates}. 
We also compare the behavior of the proposed NN with that of other NNs defined by different strategies. \cblue{In the following, we denote the interpolated VPINN as IVPINN to distinguish it from the non-interpolated VPINN \cite{kharazmi2021hp}, simply denoted as VPINN, and the standard PINN \cite{raissi2019physics}.} 

In the subsequent plots, we report by blue dots the error $\Vert u - u_H^{\cal N\!N} \Vert_{1,\Omega}$, where $u_H^{\cal N\!N}$ is the interpolated VPINN defined on the mesh $\mathcal T_H$ as in \eqref{eq:def-uHNN}, versus the size $H$ of the mesh $\mathcal T_H$.  We also show a blue solid line and a blue dashed one: the former is the regression line fitting the blue dots (possibly ignoring the first ones); its slope in the log-log plane yields the empirical convergence rate. The latter is used as a reference, since its slope corresponds  to an error decay proportional to $h^{\kt}$, \cblue{which is  the expected convergence rate of the $H^1$ error as indicated by Theorem \ref{teo:a-priori bound}, assuming that the ratio $\frac{C_h}{c_h}$ may be neglected (see Remark \ref{rem:equiv-ch-bis}).} The dashed line represents the best convergence rate we can expect from the proposed discretization scheme.

For comparison, we also report by green dots the error $\Vert u - \hat{u}^{\cal N\!N} \Vert_{1,\Omega}$, where $\hat{u}^{\cal N\!N}$ is the non-interpolated VPINN defined in Remark \ref{rem:no-interp}, and by red dots the error $\Vert u - \widetilde{u}^{\cal N\!N} \Vert_{1,\Omega}$, where $\widetilde{u}^{\cal N\!N}$ is the standard PINN proposed in \cite{raissi2019physics}, with the same architecture of the used VPINNs and the loss function computed as described in \cite{mishra2022estimates}. To obtain a fair comparison, the regularization coefficient and the ratio between the control points inside the domain and the ones on the boundary are chosen as described in \cite{mishra2022estimates}. Since we are interested in convergence rates with respect to mesh refinement, but the PINN does not require any mesh, the corresponding errors are computed by training the network with the same numbers of inputs used during the training of $u_H^{\cal N\!N}$; to be precise, the PINN is trained using the same number of collocation points as the number of interpolation nodes used by the interpolated VPINN.  

Furthermore, in order to better analyze the trade-off between the model accuracy and the training efficiency and complexity, we plot the same errors versus the dataset size. Whenever these dots, possibly after a pre-asymptotic phase, sit close to their regression line, we draw it as well in  green or red, respectively.

\medskip

\noindent \textit{Convergence test \#1: $u\in C^\infty(\bar{\Omega})$}

Consider problem \eqref{eq:model-pb} with $\Gamma_D = \{(x,y)\in \partial\Omega:x=0\text{ or }x=1\}$ and $\Gamma_N = \partial\Omega\backslash\Gamma_D$. Let us choose the following operator coefficients
\begin{equation*}
\mu(x,y) = 2+\sin(x+2y),\hspace{0.5cm} \beta(x,y) = \begin{bmatrix}
    \sqrt{x-y^2+5} \\
    \sqrt{y-x^2+5} 
\end{bmatrix}, \hspace{0.5cm} \sigma(x,y)=e^{\frac x2-\frac y3}+2\,,
\end{equation*}
and the data $f$, $g$, and $\psi$ such that the exact solution is 
\[
u(x,y) = \sin(3.2x(x-y))\cos(4.3y+x)+\sin(4.6(x+2y))\cos(2.6(y-2x)).
\]
The corresponding error decays with respect to the meshsize $H$ are shown in Fig. \ref{fig:sol5_neumann}. 

\begin{figure}[t!]
\centering 
\captionsetup{justification=centering}
\begin{subfigure}[t]{0.32\linewidth}
  \includegraphics[width=0.98\linewidth,clip]{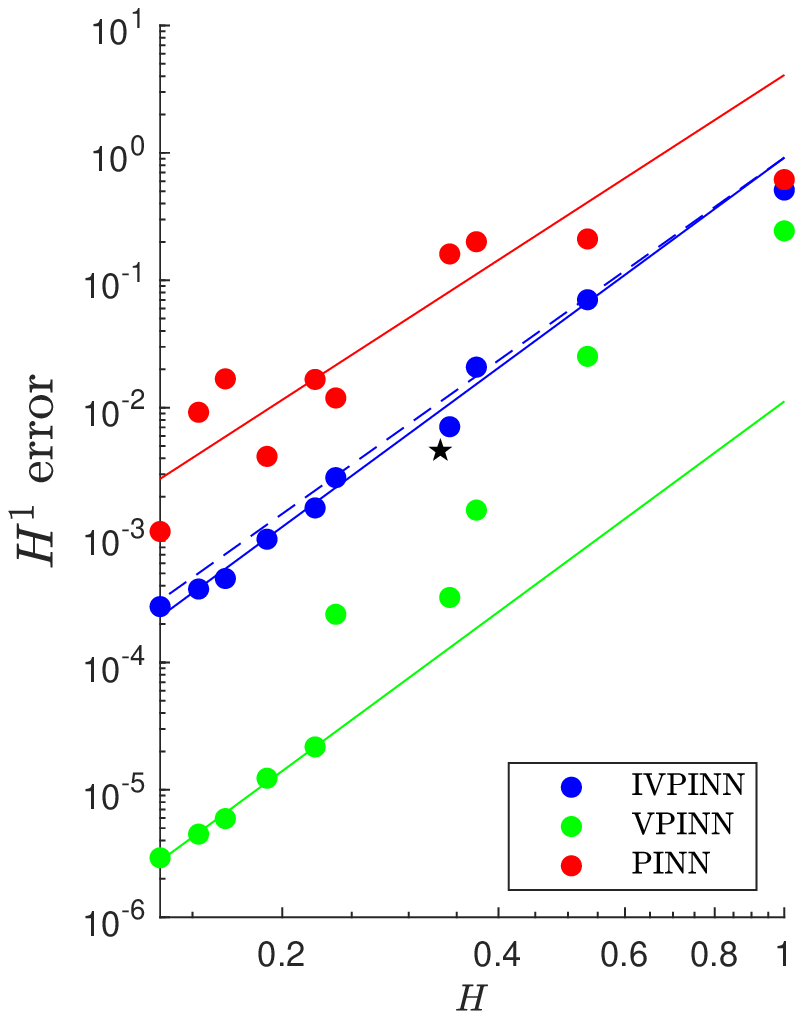} 
   \subcaption{$q=3$, $k_{\text{test}}=1$, $k_\text{int}=4$.\\ Expected decay rate: 4. \\Obtained decay rate: 4.15.}
   \label{fig:eq2_sol5n_311}
\end{subfigure}
\hfill
\begin{subfigure}[t]{0.32\linewidth}
  \includegraphics[width=0.98\linewidth,clip]{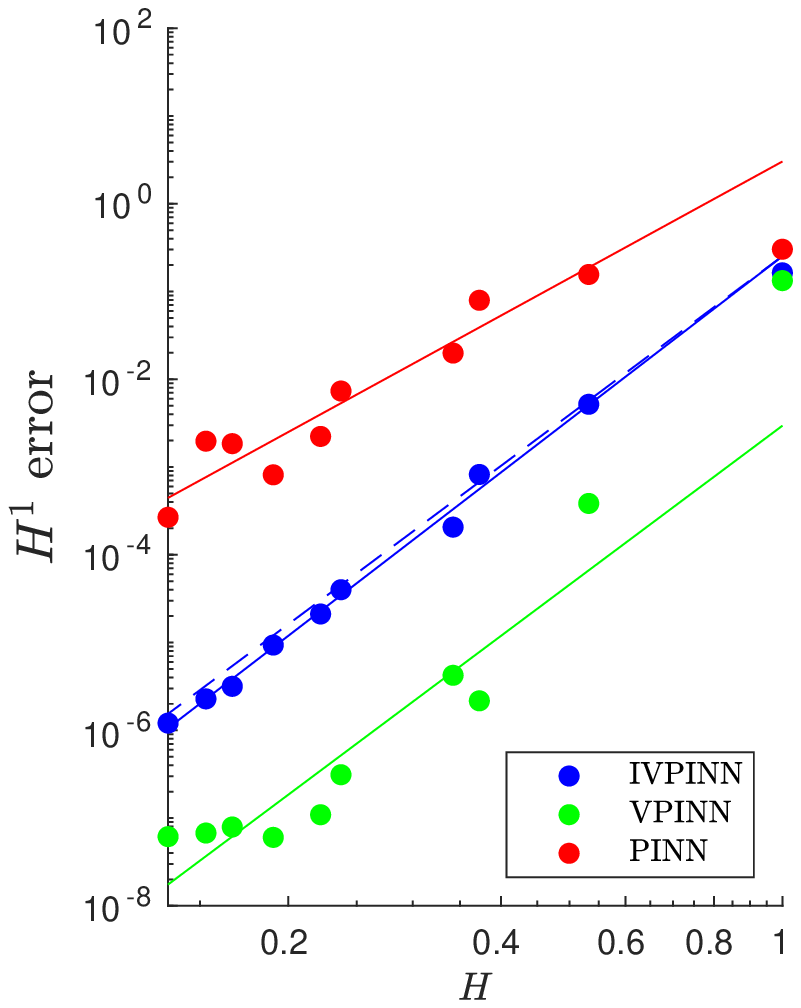} 
   \subcaption{$q=5$, $k_{\text{test}}=1$, $k_\text{int}=6$.\\ Expected decay rate: 6. \\Obtained decay rate: 6.19.}
   \label{fig:eq2_sol5n_511}
\end{subfigure}
\hfill
\begin{subfigure}[t]{0.32\linewidth}
  \includegraphics[width=0.98\linewidth,clip]{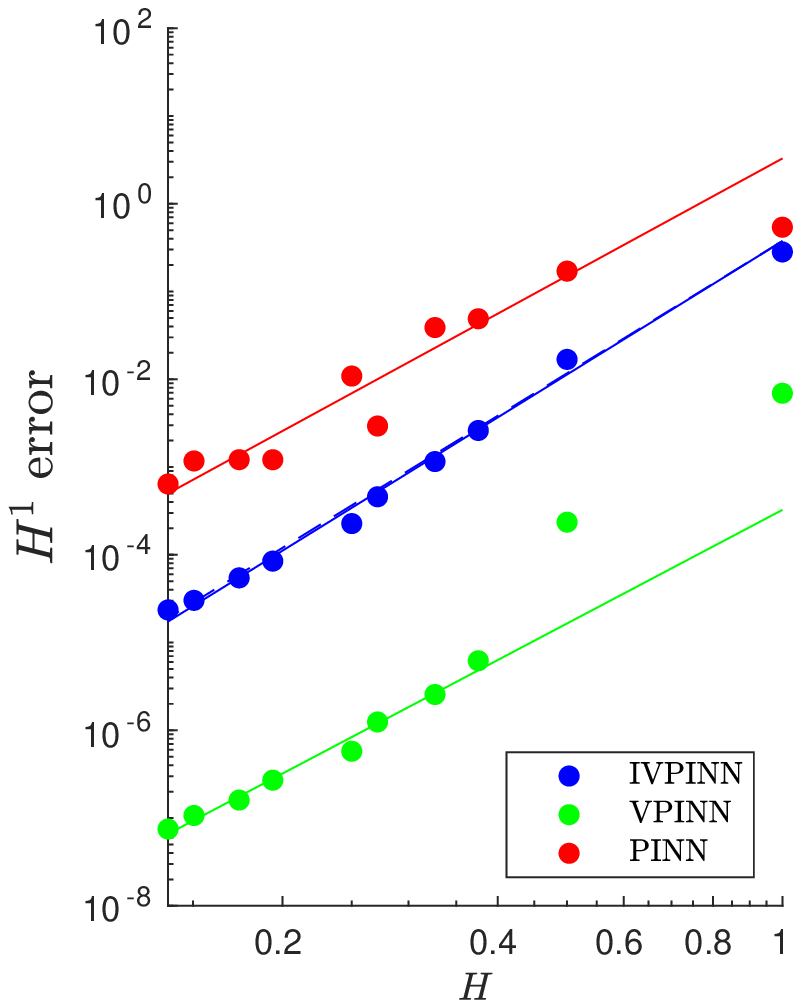}
   \subcaption{$q=5$, $k_{\text{test}}=2$, $k_\text{int}=5$.\\ Expected decay rate: 5.\\Obtained decay rate:  5.05.}
   \label{fig:eq2_sol5n_521}
\end{subfigure}
  \caption{Error decays versus $H$ for \textit{Convergence test \#1: $u\in C^\infty(\bar{\Omega})$}.}
  \label{fig:sol5_neumann}
\end{figure}
In subfigure \protect{\subref{fig:eq2_sol5n_311}}, where the \cblue{IVPINN} (blue dots) and the \cblue{VPINN} (green dots) are trained with $q=3$ and $k_{\text{test}}=1$, we observe that the points \cblue{are distributed}, possibly after an initial preasymptotic phase, \cblue{along} straight lines with slopes very close to $k_{\text{int}}=4$. \cblue{We highlight that the PINN convergence is significantly noisier and that the corresponding $H^1$ error is, on average, about 7 times the IVPINN one}. A similar phenomenon can be seen in subfigure \protect{\subref{fig:eq2_sol5n_511}}, although the finite precision of the used Tensorflow software prevents convergence to display at full for small values of $H$. In this test we use $q=5$ and $k_{\text{test}}=1$ and the regression lines for the \cblue{IVPINN} and the VPINN have slopes close to $k_{\text{int}}=6$, while the PINN accuracy is again much lower. Finally, the data in subfigure \protect{\subref{fig:eq2_sol5n_521}} are obtained with $q=5$ and $k_{\text{test}}=2$ and the \cblue{blue regression line slope is 5.05, almost coinciding with $k_{\text{int}}=5$.}

Such results highlight that, although the VPINN implementation is more complex than the PINN one, the former produces more accurate solutions than the latter, when the exact solution is regular. 

In Fig. \ref{fig:sol5_neumann_ninput}, the same error decays are expressed in terms of the number of training points, i.e., the number of neural network forward evaluations required to construct the loss function in a single epoch. Such an alternative visualization highlights that the performances of the \cblue{IVPINN} and the VPINN are very similar when trained with similar training sets. This is due to the fact that, since we stabilize the VPINN by projecting it on a space of continuous piecewise polynomials, we need fewer interpolation points (input data of the \cblue{IVPINN}) than quadrature points (input data of the VPINN) to evaluate the loss function.

\begin{figure}[t!]
\centering 
\captionsetup{justification=centering}
\begin{subfigure}[t]{0.32\linewidth}
  \includegraphics[width=0.98\linewidth,clip]{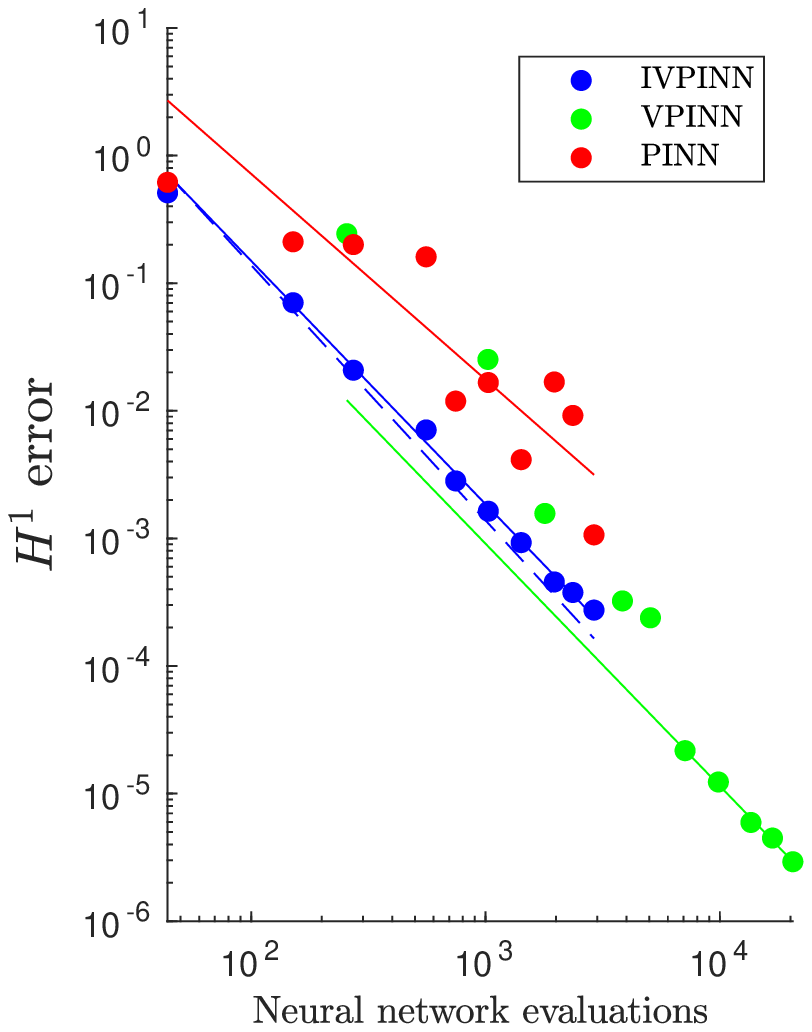} 
   \subcaption{$q=3$, $k_{\text{test}}=1$, $k_\text{int}=4$.\\ Expected decay rate: -2. \\Obtained decay rate: -1.91.}
   \label{fig:eq2_sol5n_311_ninput}
\end{subfigure}
\hfill
\begin{subfigure}[t]{0.32\linewidth}
  \includegraphics[width=0.98\linewidth,clip]{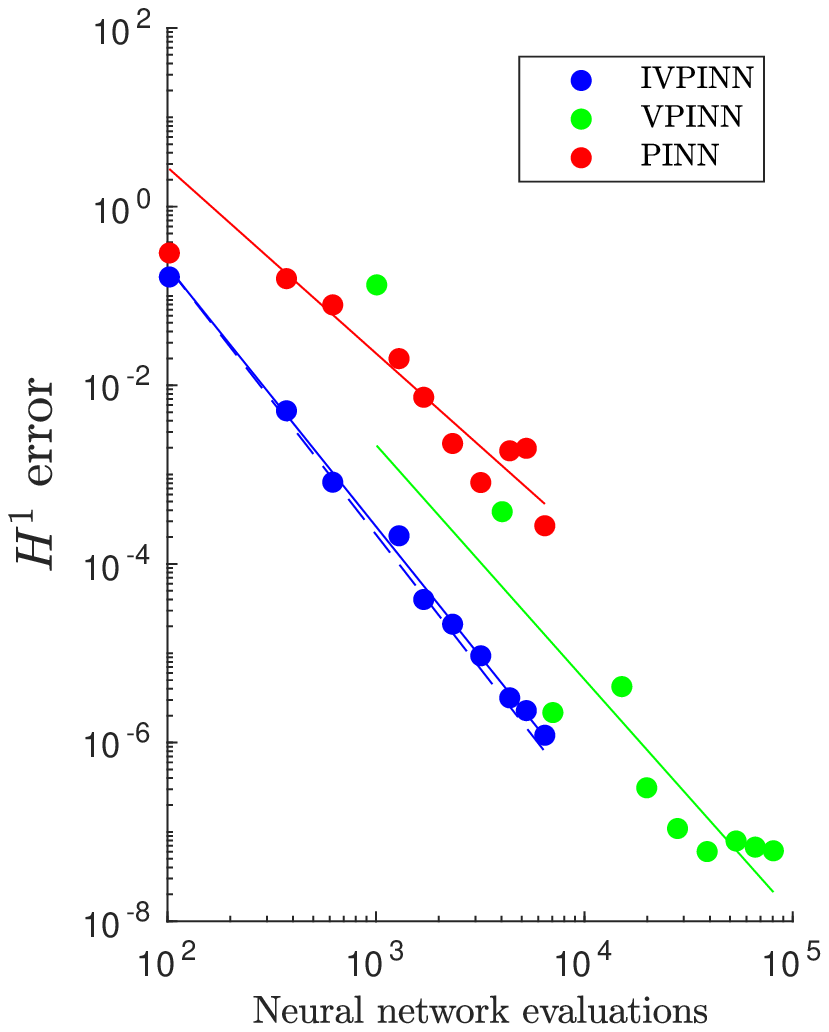} 
   \subcaption{$q=5$, $k_{\text{test}}=1$, $k_\text{int}=6$.\\ Expected decay rate: -3. \\Obtained decay rate: -2.91.}
   \label{fig:eq2_sol5n_511_ninput}
\end{subfigure}
\hfill
\begin{subfigure}[t]{0.32\linewidth}
  \includegraphics[width=0.98\linewidth,clip]{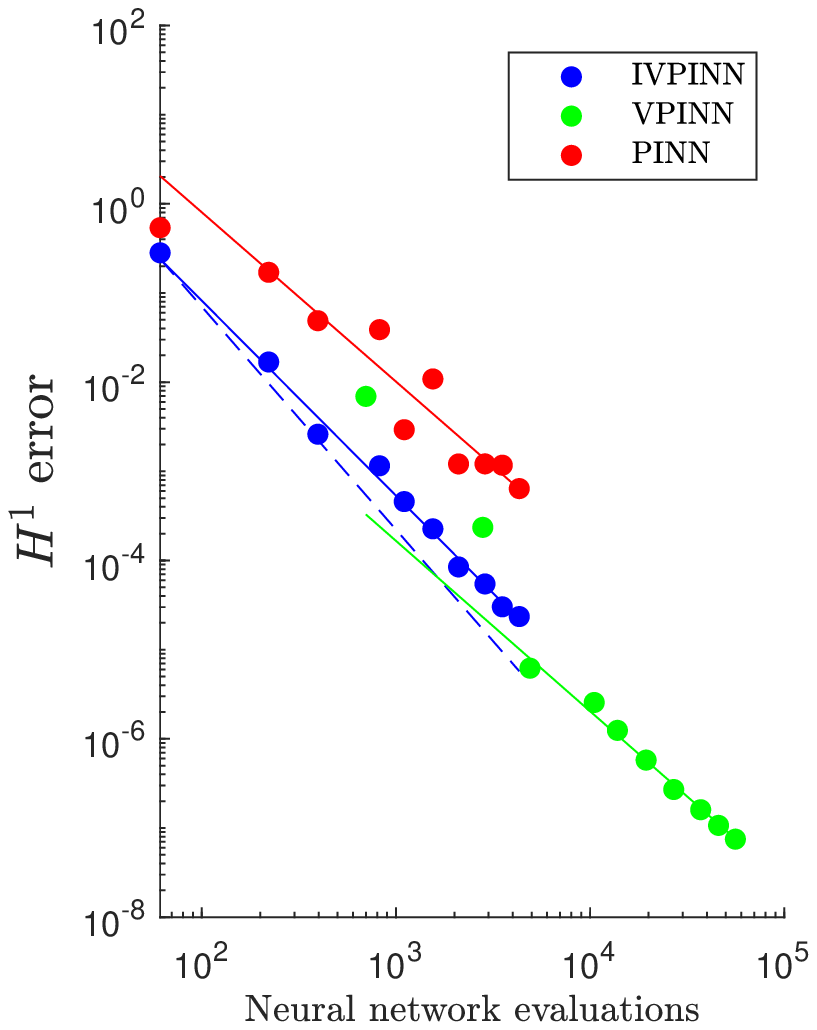}
   \subcaption{$q=5$, $k_{\text{test}}=2$, $k_\text{int}=5$.\\ Expected decay rate: -2.5.\\Obtained decay rate:  -2.19.}
   \label{fig:eq2_sol5n_521_ninput}
\end{subfigure}
  \caption{Error decays versus dataset size for \textit{Convergence test \#1: $u\in C^\infty(\bar{\Omega})$}.}
  \label{fig:sol5_neumann_ninput}
\end{figure}
}

\cblue{\begin{remark}[on the quotient $\frac{C_h}{c_h}$]\label{rem:equiv-ch-bis} 
Theorem \ref{teo:a-priori bound} indicates  that the best possible convergence rate, when the solution is regular enough, is $k_\text{int}$. However, as discussed in Remark \ref{rem:equiv-ch}, the quotient $\frac{C_h}{c_h}$ is of order $O(h^{-1})$ when test functions are picked from the Lagrange basis associated with a quasi-uniform triangulation and the weights $\gamma_i$ are equal to 1. In this case, the term $\left(1+\frac{C_h}{c_h}\right)$ in \eqref{eq:a-priori bound} reduces the predicted convergence by exactly one order.

On the other hand, in Fig. \ref{fig:sol5_neumann}, we have shown cases where the order of convergence is optimal. Such a behavior is related to the fact that the loss $R_h(u^{\cal N\!N})$ decays much faster than expected, in the considered cases, namely at least as $O(h^8)$. Therefore, when $h$ is small enough, the term $C_h R_h(u^{\cal N\!N})$ in Lemma \ref{lem:intermedio1} can be neglected, and the predicted convergence rate is not affected by the presence of the quotient $\frac{C_h}{c_h}$.
\end{remark}}

\medskip

\noindent \cblue{\textit{Convergence test \#2: $u\in H^{5/3-\varepsilon}(\Omega)$}

Let us now focus on a less smooth solution, whose regularity is commonly found in domains with reentrant corners. The problem is characterized by $\Gamma_D=\partial \Omega$, $\mu=1$, $\beta=[2,3]^T$, $\sigma=4$, whereas the forcing term and boundary conditions are such that the exact solution is, in polar coordinates,
\[
u(r,\theta) = r^\frac 23 \sin\left(\frac 23\left(\theta+\frac \pi2\right)\right).
\]
Since $u\in H^{5/3-\varepsilon}(\Omega)$ for any $\varepsilon>0$, we expect a convergence rate close to $2/3$; indeed, $k_{\text{int}}\ge 4$ and the rate of convergence is always limited by the solution regularity as expected. The error decays are shown in Fig. \ref{fig:sol4}. Notice that the \cblue{IVPINN} is even more stable and accurate than the VPINN trained on the same mesh (i.e., with more input data). The PINN behaves better in this test case than in the previous one, but still the accuracy is worse than the one provided by our \cblue{IVPINN}.
\begin{figure}[t!]
\centering 
\captionsetup{justification=centering}
\begin{subfigure}[t]{0.32\linewidth}
  \includegraphics[width=0.98\linewidth,clip]{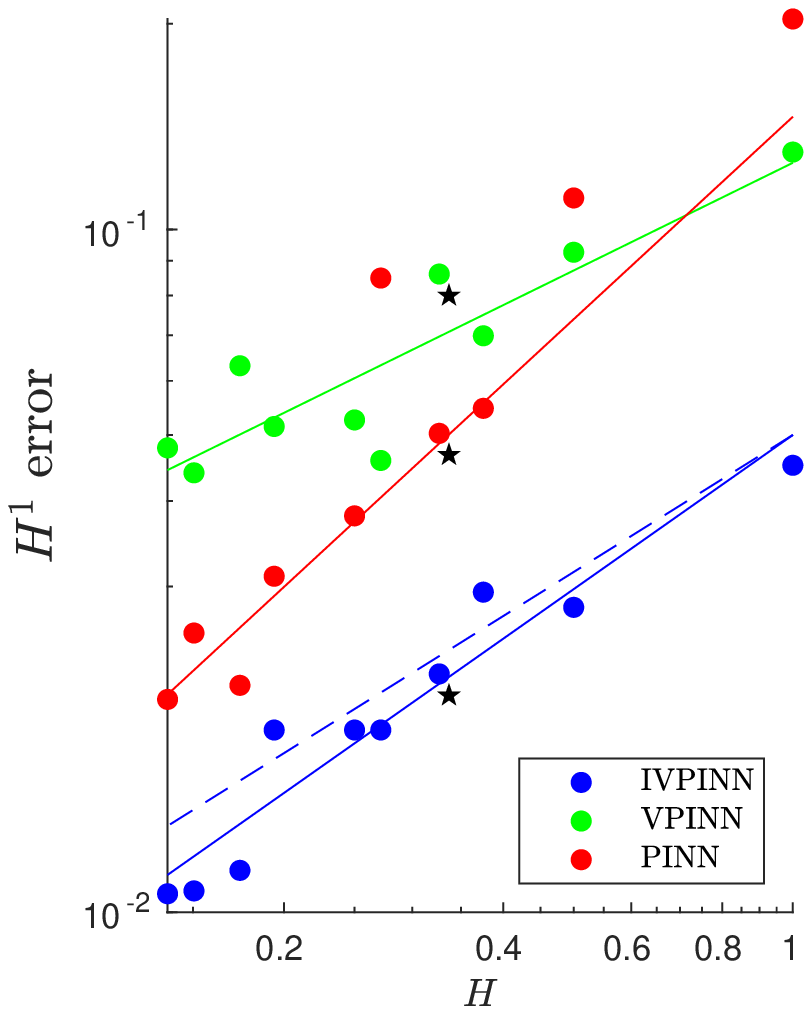} 
   \subcaption{$q=3$, $k_{\text{test}}=1$, $k_\text{int}=4$.\\ Expected decay rate: 2/3. \\Obtained decay rate: 0.75.}
   \label{fig:eq1_sol4_311}
\end{subfigure}
\hfill
\begin{subfigure}[t]{0.32\linewidth}
  \includegraphics[width=0.98\linewidth,clip]{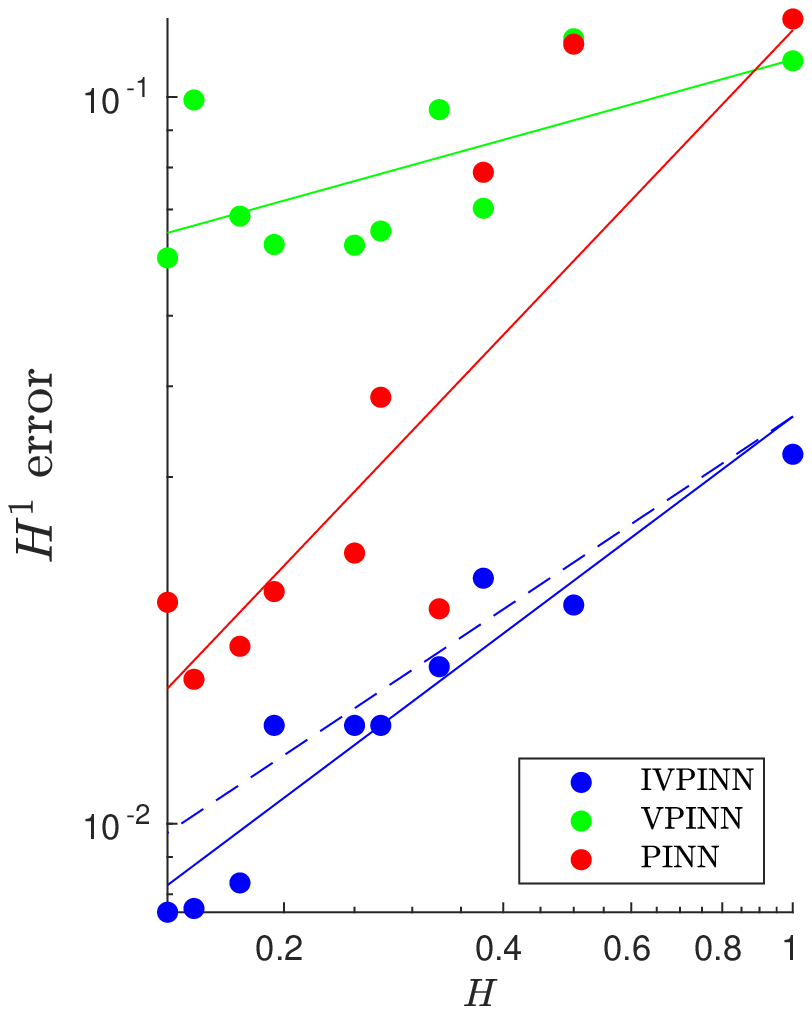} 
   \subcaption{$q=5$, $k_{\text{test}}=1$, $k_\text{int}=6$.\\ Expected decay rate: 2/3. \\Obtained decay rate: 0.75.}
   \label{fig:eq1_sol4_511}
\end{subfigure}
\hfill
\begin{subfigure}[t]{0.32\linewidth}
  \includegraphics[width=0.98\linewidth,clip]{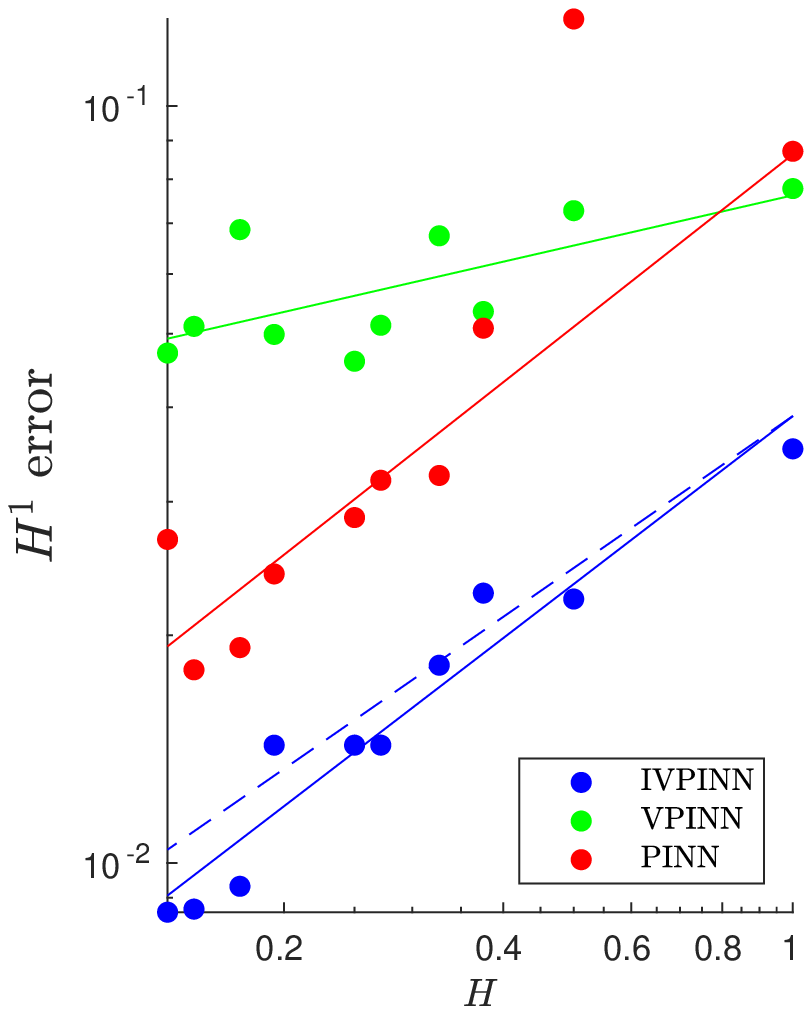}
   \subcaption{$q=5$, $k_{\text{test}}=2$, $k_\text{int}=5$.\\ Expected decay rate: 2/3.\\Obtained decay rate:  0.74.}
   \label{fig:eq1_sol4_521}
\end{subfigure}
  \caption{Error decays versus $H$ for \textit{Convergence test \#2: $u\in H^{5/3-\varepsilon}(\Omega)$}.}
  \label{fig:sol4}
\end{figure}

It is also interesting to analyze the behavior of the loss function and of the error during training as documented in Fig. \ref{fig:sol4_performances}, where the first 3000 epochs are performed with the ADAM optimizer, while the remaining ones with the BFGS optimizer.} \cblue{Such plots correspond to the loss function and the $H^1$ error associated with the dots marked by the black stars in Fig. \protect{\subref{fig:eq1_sol4_311}}.} \cblue{It can be noted that the \cblue{IVPINN} and the VPINN initially converge very fast with the ADAM optimizer; eventually, after the initial phase in which both the loss and the error decrease, the error reaches a constant value despite the loss function keeps diminishing. This implies that there exist other sources of error that prevail when the loss function decays. On the other hand, using a standard PINN, one observes that the convergence of the loss and the error is much slower than for the VPINNs, and the second-order optimizer is needed to converge to an accurate solution. The average epoch execution time is approximately 0.0599 seconds for the PINN, 0.0587 seconds for the VPINN, and  0.0479 seconds for the \cblue{IVPINN}. Such a gain is due to the fact that the model derivatives are computed via automatic differentiation in the non-interpolated models, while the gradient of the \cblue{IVPINN} can be computed by a simple matrix-vector multiplication. Note that the gain increases when higher derivatives are involved in the PDE.
\begin{figure}[t!]
\centering 
\captionsetup{justification=centering}
\begin{subfigure}[t]{0.49\linewidth}
  \includegraphics[width=0.98\linewidth,clip]{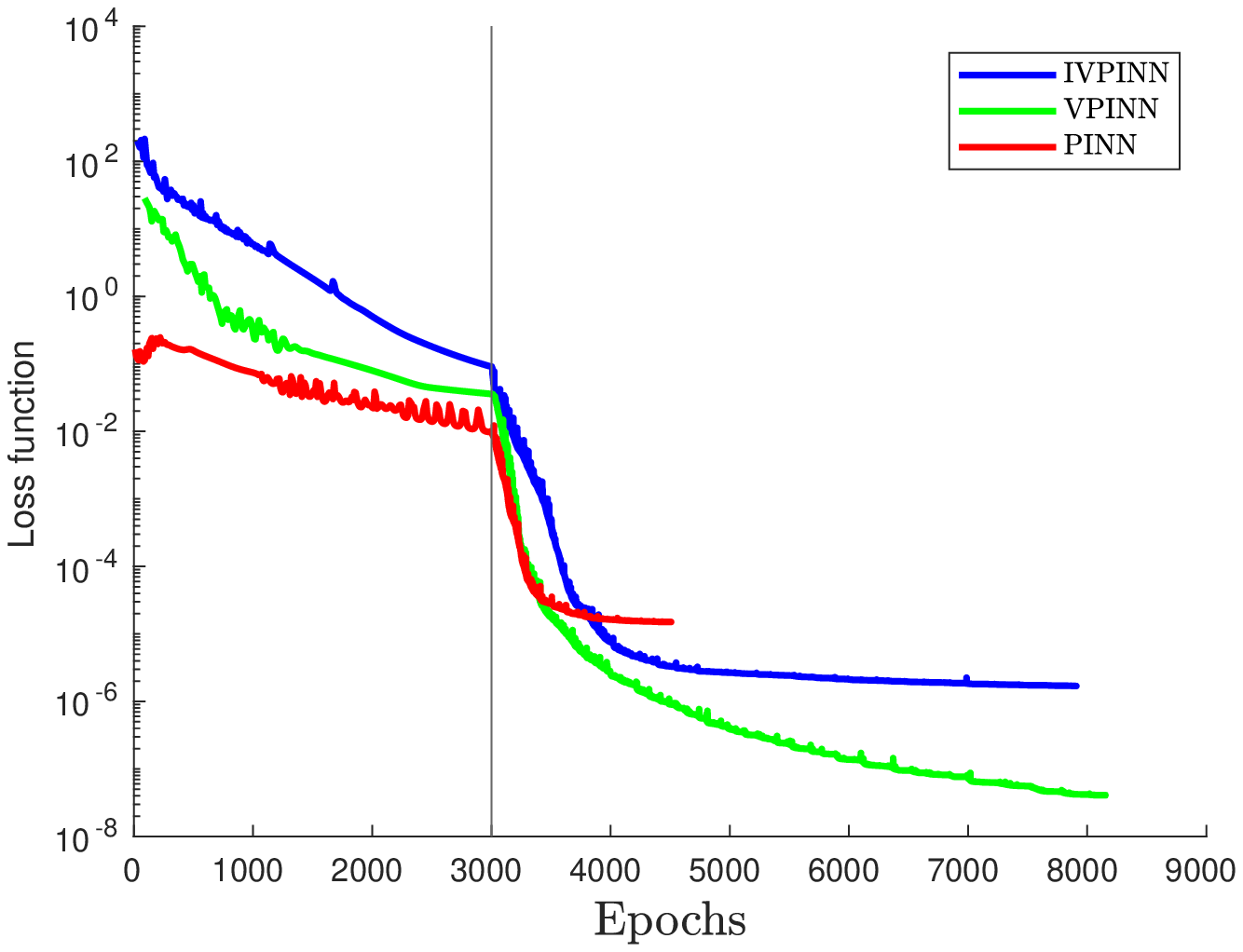} 
   \label{fig:loss_epochs}
\end{subfigure}
\hfill
\begin{subfigure}[t]{0.49\linewidth}
  \includegraphics[width=0.98\linewidth,clip]{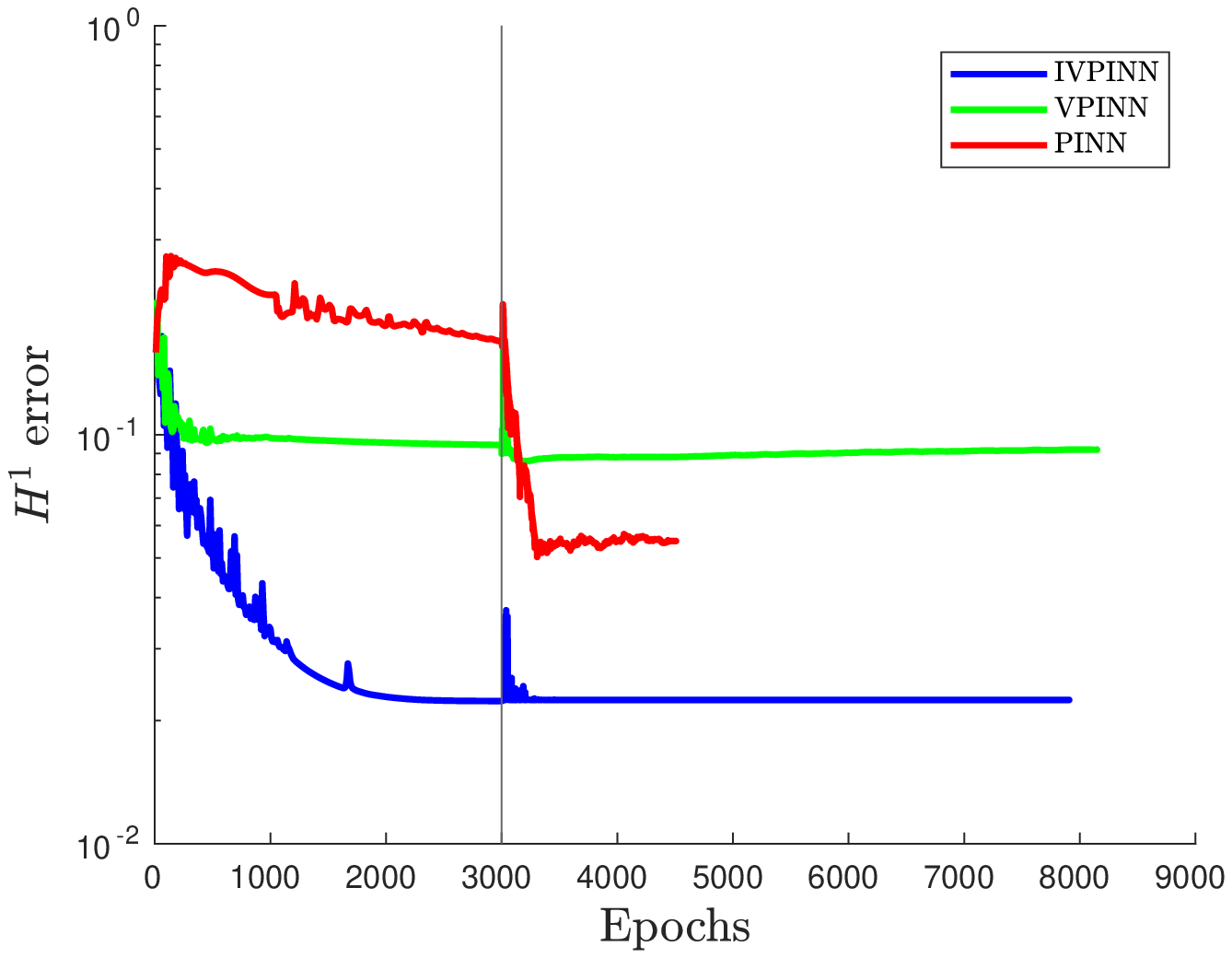} 
   \label{fig:errors_epochs}
\end{subfigure}

\begin{subfigure}[t]{0.49\linewidth}
  \includegraphics[width=0.98\linewidth,clip]{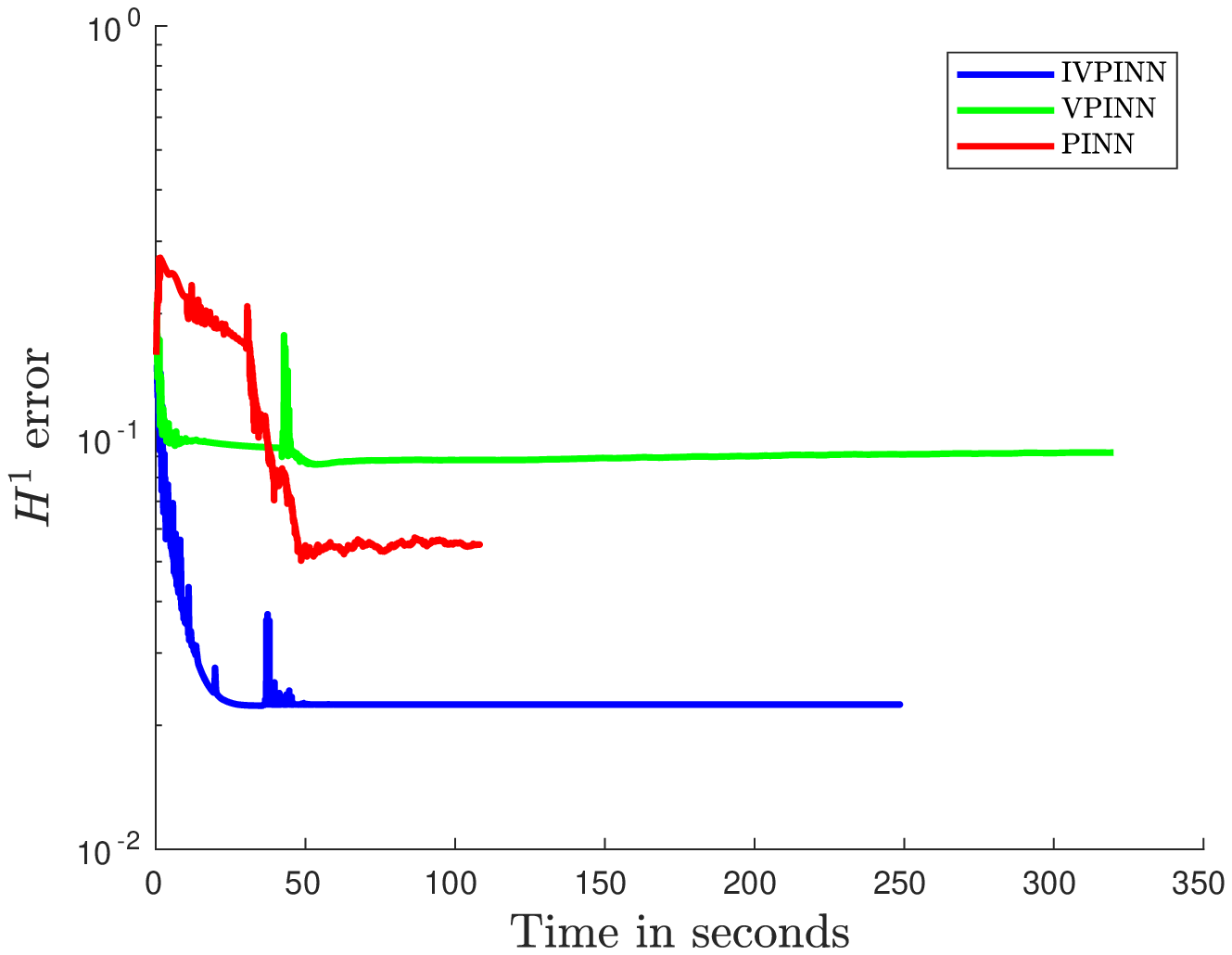} 
   \label{fig:errors_time}
\end{subfigure}
  \caption{Top row: loss function (left) and $H^1$ error (right) evaluations as functions of the number of executed epochs. The first 3000 epochs are performed with the ADAM optimizer, the subsequent ones with the BFGS optimizer. \\ Bottom row:  $H^1$ error as a function of the elapsed time.}
  \label{fig:sol4_performances}
\end{figure}
}
\cblue{\subsection{How the VPINN dimension affects accuracy} \label{sect:error_vs_dimension}
We now focus on the dependence of the error on the neural network dimension. For the sake of simplicity, we fix the problem discretization and vary only the number of layers and the number of neurons in each layer, assuming that each layer contains exactly the same number of neurons. The considered domain, parameters, forcing term and boundary conditions are the ones described in \textit{Convergence test \#1}. The VPINN is trained with piecewise linear test functions ($k_{\text{test}}=1$) and quadrature rules of order $q=3$ \cblue{on the finest mesh used to produce Fig. \protect{\subref{fig:eq2_sol5n_311}} and on the mesh associated with the blue dot close to the black star in the same figure.}

We can observe, in Fig. \ref{fig:error-vs-dimension}, that the error is very high for small networks, but then it rapidly decreases while increasing the number of neurons in each layer, until a plateau is reached depending on the chosen problem discretization. Essentially, on both meshes, 3 layers with 10 neurons each suffice to achieve the lowest possible discretization error for the given loss function.

This analysis confirms that the error decays reported in Sect.  \ref{sect:error_vs_h} are all insensitive to the neural network hyper-parameters, as they have been obtained by a large neural network (5 layers with 50 neurons each).
Such results validate the assumption made in Section \ref{sec:error_estimates} about  the neural network, namely that its dimension -- provided it is sufficiently large -- does not influence the predicted convergence rate. 

\begin{figure}[t!]
\centering 
\captionsetup{justification=centering}

\begin{subfigure}[t]{0.49\linewidth}
  \includegraphics[width=0.98\linewidth,clip]{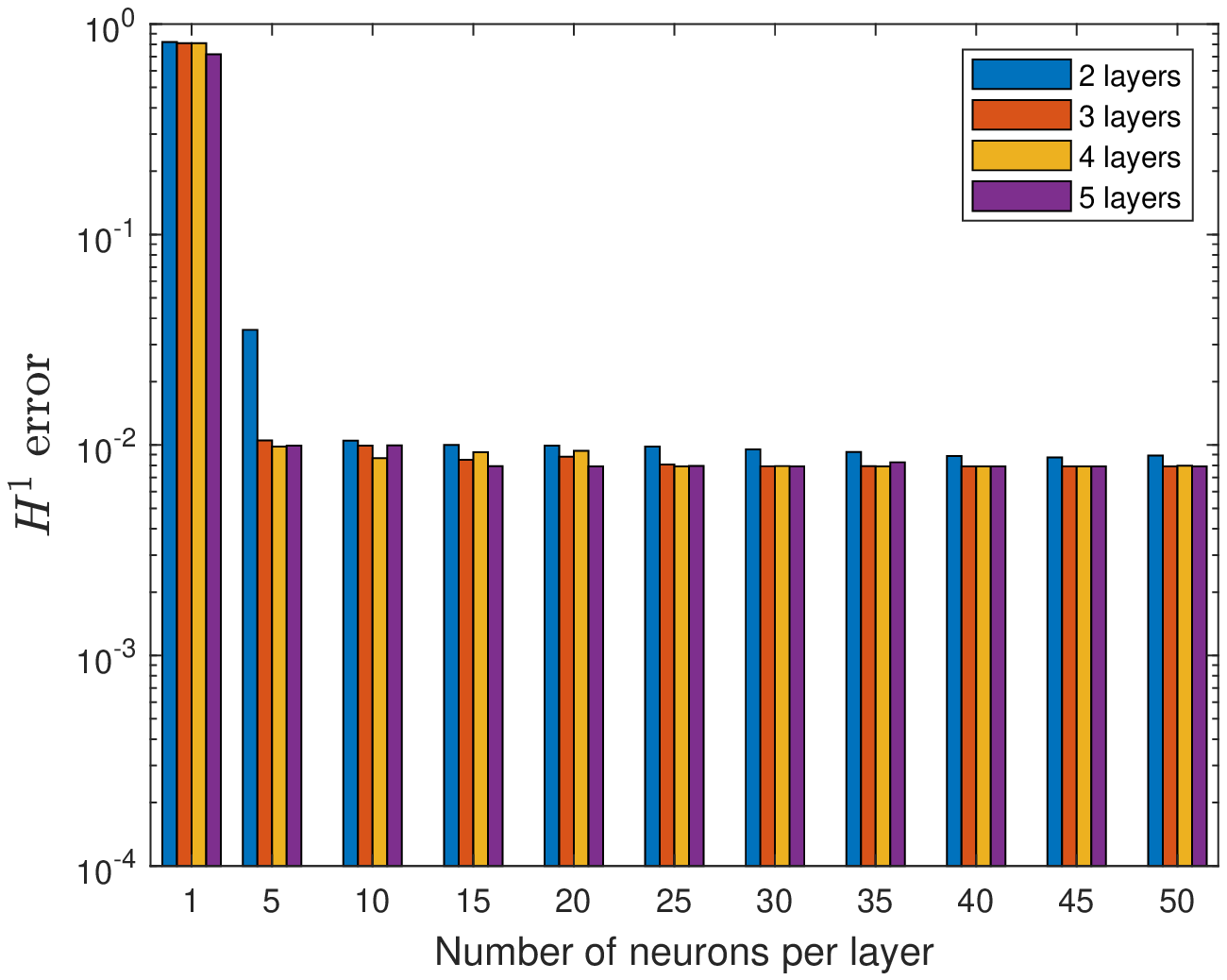} 
 \subcaption{VPINN trained on the mesh used to compute \cblue{the blue point close to the black star} in Fig. \protect{\subref{fig:eq2_sol5n_311}}.}
   \label{fig:error_vs_dimensions_interm} 
\end{subfigure}
\hfill
\begin{subfigure}[t]{0.49\linewidth}
  \includegraphics[width=0.98\linewidth,clip]{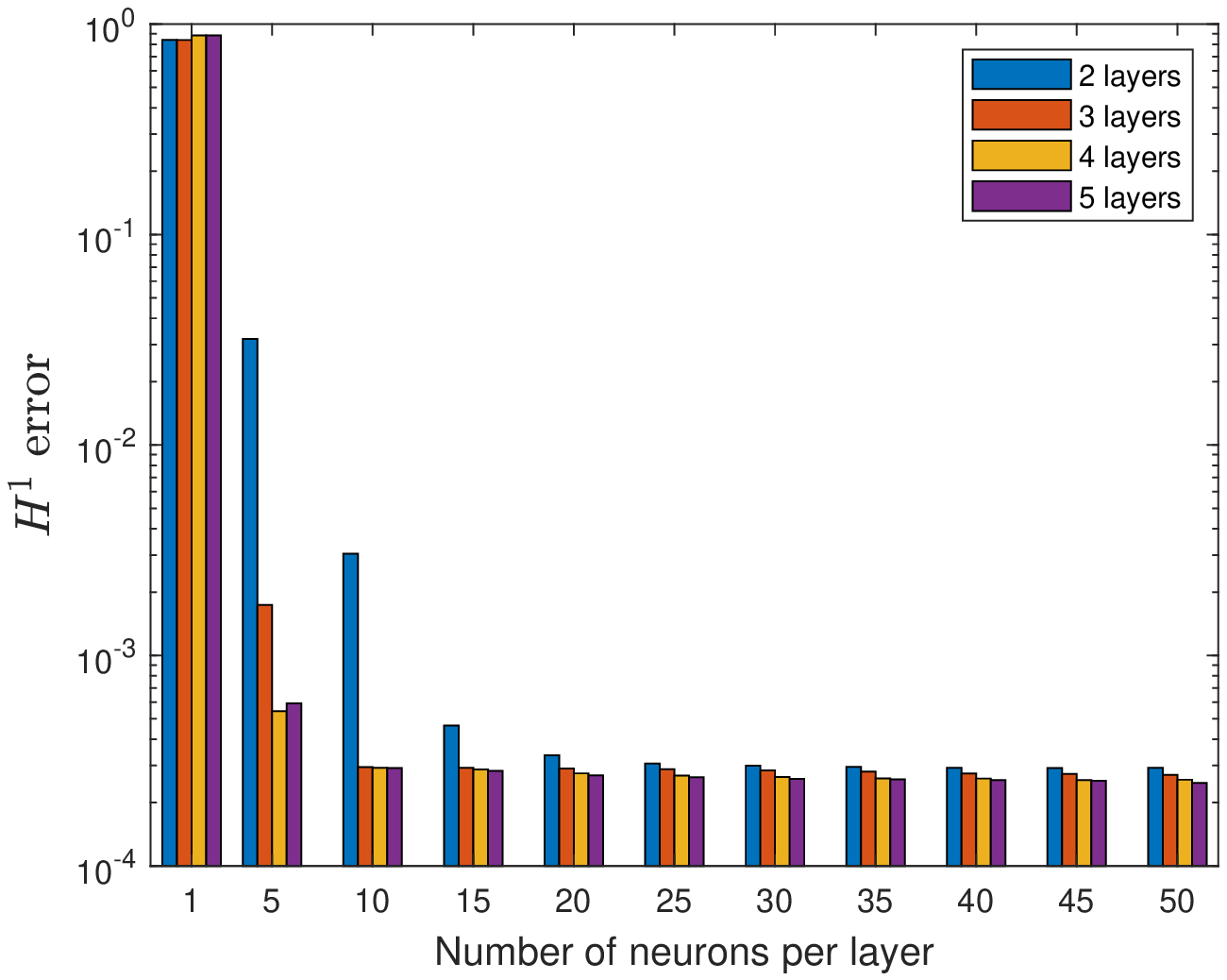} 
   \subcaption{VPINN trained on the finest mesh used in Fig. \protect{\subref{fig:eq2_sol5n_311}}.}
   \label{fig:error_vs_dimensions_fine} 
\end{subfigure}
  \caption{$H^1$ error varying the number of layers and the number of neurons in each layer of the neural network.}
  \label{fig:error-vs-dimension}
\end{figure}

}

\subsection{On the importance of the inf-sup condition} \label{sect:inf-sup-section}
In this section we show that the inf-sup condition, assumed in Proposition \ref{cor:discrete-inf-sup} to derive the a priori error estimate, is crucial in order to avoid spurious modes in the numerical solution. To prove such a claim, let us consider the simplest one-dimensional Poisson's problem with zero forcing term and zero Dirichlet boundary conditions:
\begin{equation}\label{eq:sol0}
\left\{
\begin{aligned}
- &u''  = 0 && \text{in \ } \Omega\,, \\
& \ \, u=0 &&  \text{on \ } \Gamma_D \,,
\end{aligned}\right.
\end{equation}
where $\Omega=(0,1)$ and $\Gamma_D=\{0,1\}$. For the sake of simplicity, we use piecewise linear test functions and quadrature rules of order $q=3$.
Note that since both $U^{\cal N\!N}$ and $U_H$ contain the exact solution $u\equiv 0$, it is always possible to obtain a numerical solution that is identical to the exact one (up to numerical precision). 

\begin{figure}[t!]
\centering 
\captionsetup{justification=centering}
\begin{subfigure}[t]{0.49\linewidth}
  \includegraphics[width=0.98\linewidth,clip]{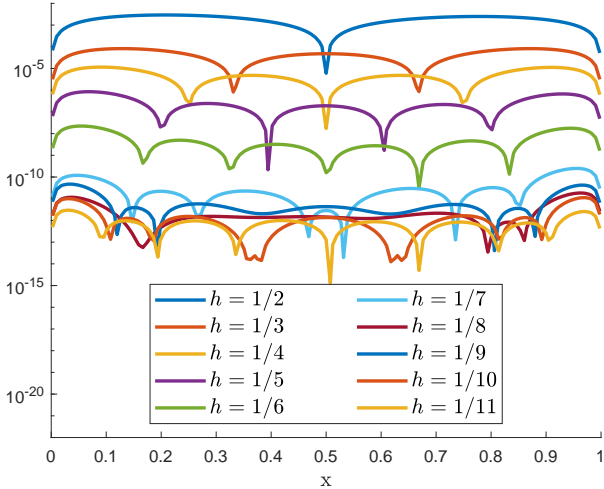} 
   \label{fig:sol0_noinfsup}
\end{subfigure}
\hfill
\begin{subfigure}[t]{0.49\linewidth}
  \includegraphics[width=0.98\linewidth,clip]{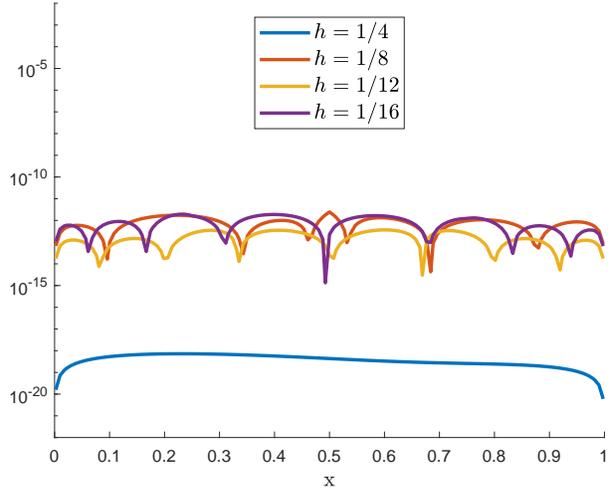}
   \label{fig:sol0_infsup}
\end{subfigure}
  \caption{Numerical solutions of problem \eqref{eq:sol0} computed with different meshes, solutions obtained with non-interpolated neural networks (left) and with interpolated neural networks (right). Quadrature rule order $q=3$. Test functions order $k_{\text{test}}=1$. }
  \label{fig:sol0}
\end{figure}

Let us denote by $u_\delta$ any discrete solution defined in Sect. \ref{sec:sub_discretization}, namely, either a solution $u_H^{\cal N\!N}$ obtained by interpolated VPINNs, or a solution $\hat{u}^{\cal N\!N}$ obtained by non-interpolated VPINNs. 
These discrete solutions are represented in Fig. \ref{fig:sol0} in logarithmic scale to allow a direct comparison. In order to avoid numerical issues due to the logarithmic scale of the plot when $u_{\delta}$ gets close to 0, a truncation procedure is applied.

The functions ${u}_\delta$ produced by non-interpolated networks are represented in the left plot of Fig. \ref{fig:sol0}. Each one is obtained by minimizing the loss function up to machine accuracy; despite this, when the mesh is fairly coarse the discrete solution is significantly different from the null solution. Indeed, the initial weights in the training process are non-zero, and the minimization process is under-determined, thereby allowing the existence of non-zero global minima. Refining the mesh, the approximation improves up to a maximum precision imposed by the chosen network architecture and the Tensorflow deep learning framework.

Conversely, the plots in the subfigure on the right-hand side of Fig. \ref{fig:sol0}, produced by interpolated networks, clearly indicate that the obtained discrete solutions are numerically zero, irrespective of the meshsize. Note that the case $h=1/4$ differs from the others since here the interpolation mesh ${\cal T}_H$ is formed by just one element. The corresponding function $u_H^{\cal N\!N}$ is thus differentiable everywhere and the used gradient-based optimizers are able to minimize the loss function more effectively. We highlight that, since $q=3$ and $k_{\text{test}}=1$, we have to choose $H=4h$ to satisfy the inf-sup condition.

\cblue{To illustrate the mechanism that may lead to the onset of spurious modes, in Appendix \ref{sect:stability_large_dataset}   we provide an analytical example of a neural network which significantly differs from a PDE solution, yet it is a global minimizer of the corresponding loss function.  }

These results, although obtained in overly simple functional settings, show the potential existence of uncontrolled components in the discrete solutions obtained by non-interpolated neural networks.
 In more complex scenarios, the presence of spurious modes may be even more pronounced.
 \cblue{ In practice, as observed in Fig. \ref{fig:sol5_neumann}, when the PDE solution is smooth enough non-interpolated solutions appear to be more accurate than the corresponding solutions obtained by interpolated neural networks using the same test functions}; however, a rigorous analysis of NN-based discretization schemes should also cope with the presence of spurious components, which we have avoided by resorting to an inf-sup condition.  We believe that these observations shed new light on the use of deep learning in numerical PDEs.

\cblue{\section{Application to nonlinear parametric problems} \label{sect:parametric}
In the previous sections, we investigated the features of the proposed VPINN discretization for the linear boundary-value problem \eqref{eq:model-pb}. Hereafter, we provide an application where the nonlinear nature of neural networks can be exploited at best. It is well known that solving nonlinear PDEs by PINNs or VPINNs comes at little extra cost with respect to linear PDEs, since nonlinearities just impact the computation of the loss function. Similarly, parametric problems can be easily and efficiently handled by neural networks, even when the dependence of the solution upon the parameters is nonlinear. Indeed, it is enough to add as many inputs as the number of parameters in the definition of the network, and train it on a proper subset of the parameter space.

To illustrate the behavior of our VPINN in these situations, let us consider the following nonlinear parametric equation:
\begin{equation}\label{eq:model-pb-parametric}
\begin{cases}
-\nabla \cdot (\mu \nabla u) + \boldsymbol{\beta}\cdot \nabla u + \sigma\,  {\rm e}^{-pu^2} =f & \text{in \ } \Omega\,, \\
u=g & \text{on \ } \partial\Omega \,, 
\end{cases}
\end{equation}
where $\mu, \beta, \sigma, f$ and $g$ are suitably smooth functions, and $p\in\Omega_p\subset \mathbb R$ is an additional parameter. 
Our goal is to train a neural network to compute the numerical solution for any given value of $p$ \cblue{in a prescribed parametric domain $\Omega_p$.} 

We fix $\Omega=(0,1)^2$, $\mu=1$, $\beta=[2,3]^T$, $\sigma=4$ and choose $f=f(\cdot,p)$ and $g= g(\cdot,p)$ such that the exact solution is
\[
u(x,y;p) = \frac{\cos\left(5\left(px+\frac y2\right)\right)}{1+p} + \left(x+\frac y2\right)^2.
\] 
To approximate such a solution in the parametric domain $\Omega_p=[0.5,2]$, we consider a neural network with three inputs ($x,y,p$) and we train it using the loss function:
\begin{equation}\label{eq:parametric-loss-function}
R_h^2(w) = \sum_{p\in\Omega_p^\#}\sum_{i \in I_h} r_{h,i;p}^2(w) \, \gamma_i^{-1}\,,
\end{equation}
where $\Omega_p^\#=\{p_1,...,p_{n_p}\}\subset\Omega_p$ is a finite set of parameter values, and the residuals $r_{h,i;p}$ are defined as in \eqref{eq:residuals}, considering the new equation. In this numerical test $\Omega_p^\#$ contains $n_p=13$ equally spaced values $0.5=p_1<p_2<...<p_{n_p}=2$. 

After the training phase, the neural network can be evaluated at new parameter values to analyze its accuracy. The error diagram is presented in Fig. \ref{fig:parametric-error}: the blue line is computed as the $H^1$ error between the exact solution $u=u(\cdot ,p)$ and the corresponding numerical solution, while the red dots represent the error associated with parameters in $\Omega_p^\#$. Despite the small number of parameter values used during training, the model provides accurate solutions for the whole range of parameter values. Note that the error increases for larger values of $p$ because the solution is more and more oscillating as $p$ increases.

\begin{figure}[t!]
\centering 
\captionsetup{justification=centering}
  \includegraphics[width=0.6\linewidth]{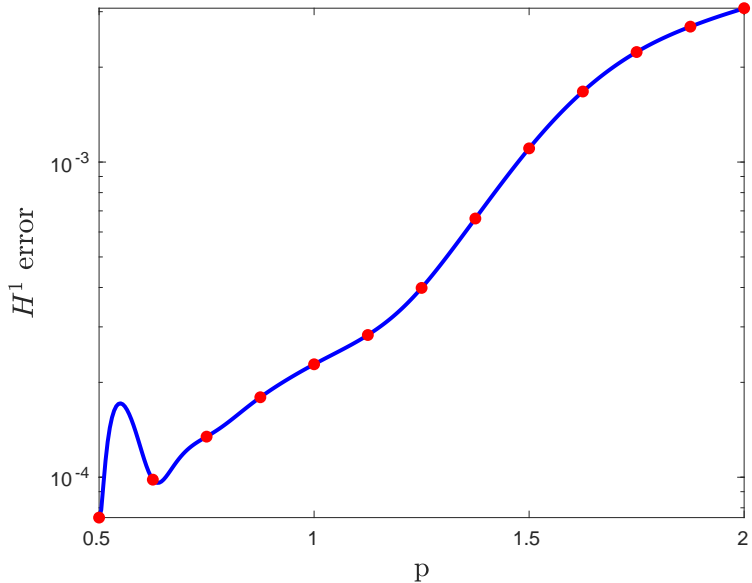}
  \caption{$H^1$ error for different values of the parameter $p$. The VPINNs are trained on $n_p=13$ parameters. The red dots are associated with parameter values used in the training phase}
  \label{fig:parametric-error}
\end{figure}
}

%% file: conclusions.tex
\section{Conclusions}\label{sec:conclusions}

We have investigated VPINN methods for elliptic boundary-value problems, in what concerns the choice of test functions and quadrature rules. The aim was the derivation of rigorous a priori error estimates for some projection of the neural network solution. The neural network is trained using as test functions finite-element nodal Lagrange basis functions of degree $k_\text{test}$ on a mesh ${\cal T}_h$, where  Gaussian quadrature rules of order $q$ are applied in each element of ${\cal T}_h$. For a fixed neural network architecture with tanh activation function, we studied how the error in the energy norm depends upon the mesh parameter $h$, for different values of $k_\text{test}$ and $q$.

Error control was obtained for the finite-element interpolant of degree $\kt=q+2-k_\text{test}$ of the neural network on an auxiliary mesh ${\cal T}_H$; such an interpolation enters also in the definition of the residuals which are minimized through the loss function. A key ingredient in the error control is the validity of an inf-sup condition between the spaces of test functions and interpolating functions. Indeed, the neural network solution might be affected by spurious modes due to the under-determined nature of the minimization problem, as we documented for a problem with zero data; instead, the onset of such modes is prevented by the adopted interpolation procedure.

Our analysis reveals that the convergence rate in the energy norm is at least of order $q+1-k_\text{test}$ for sufficiently smooth functions, and it increases to $q+2-k_\text{test}$ when the value of the loss function obtained by minimization is sufficiently small. The main message stemming from the analysis is that it is convenient to choose test functions of the lowest degree $k_\text{test}=1$ in order to get the highest convergence rate for a fixed quadrature rule. Furthermore, for smooth solutions the convergence rate may be arbitrarily increased by increasing the precision of the quadrature rule, although the realization of this theoretical statement is hampered in practice by the finite precision of machine arithmetics.

\cblue{We also investigated the influence of the neural network hyperparameters on the overall accuracy of the discretization, and we found that a small network with few layers and neurons suffices to reach accuracies of practical interest. To stay on the safe side, we used a larger network in our experiments, thereby obtaining results that are essentially independent of the network hyperparameters.

For the sake of comparison, we also implemented a standard VPINN without projection upon piecewise polynomials, as well as a standard PINN trained with the same number of inputs as those used in training our VPINN. Interestingly, we experimentally observed that in general the error decay rate for the non-interpolated neural network solution replicates the one theoretically predicted for the interpolated network. The PINN solutions appear to be less accurate and noisier than the interpolated VPINN's.

We have shown that interpolated VPINNs are able to efficiently solve nonlinear parametric problems without the need for additional nonlinear solvers or globalization methods, due to their intrinsic nonlinear nature. The VPINN can be trained in an off-line phase on a subset of the parameter domain and then efficiently evaluated on-line on any other parameter value. This is a key difference between the proposed method and standard numerical techniques such as FEM, even if the solution is sought in the same finite dimensional space. Indeed,  the latter would require some iterative technique to handle nonlinearities, as well as some form of interpolation/extrapolation to get the solution for the whole range of parameters. All this is provided for free by the NN machinery.

Possible extensions of this work are related to the investigation of more advanced neural networks architectures to improve the method accuracy and efficiency \cite{rodriguez2021physics,gao2021phygeonet} or to more complex problems. Indeed, neural networks are known to be able to manage very high-dimensional problems, overcoming the so-called curse of dimensionality, therefore we expect them to be able to efficiently solve parametric PDEs with multiple parameters \cite{gao2021phygeonet} or high-dimensional PDEs \cite{han2018solving,LanthalerMishraKarniadakis2021}. Other possible applications are related, for instance, to inverse problems \cite{chen2020electromagnetic} or integration between PDEs and data \cite{chen2021physics}.}

\bigskip
\noindent
{\bf Declarations.}  The authors performed this research in the framework of the Italian MIUR Award ``Dipartimenti di Eccellenza 2018-2022" granted to the Department of Mathematical Sciences, Politecnico di Torino (CUP: E11G18000350001). The research leading to this paper has also been partially supported by the SmartData@PoliTO center for Big Data and Machine Learning technologies.
SB  was supported by the Italian MIUR PRIN Project 201744KLJL-004, CC was supported by the Italian MIUR PRIN Project 201752HKH8-003. 
The authors are members of the Italian INdAM-GNCS research group.

The authors have no relevant financial or non-financial interests to disclose.

The datasets generated during and/or analysed during the current study are available from the corresponding author on reasonable request.

%% file: interpolation_construction.tex
\cblue{\section*{Appendix}
\section{Construction of the interpolation operator}\label{sec:Ih_construction}
In this section we provide details on the practical construction of the operator ${\cal I}_H: \Cr^0(\bar{\Omega}) \to U_H$ introduced in Section \ref{sec:sub_discretization}.

Since $U_H$ is the linear subspace of $U$ containing all the piecewise polynomial of degree $k_{\text{int}}$ defined over ${\cal T}_H$, there exists a Lagrange basis $\{\hat \varphi_i\}_{i=1}^{n_I}$ such that $U_H=\text{span}\{\hat \varphi_i:i=1,...,n_I\}$, which is associated with a corresponding set of points $\{\mathbf x_i\}_{i=1}^{n_I}\subset \Omega$. The basis functions satisfy the relations $\hat\varphi_i(\mathbf x_j)=\delta_{i,j},\forall i,j=1,...,n_I$. Therefore, the operator ${\cal I}_H$ maps the generic function $v\in\Cr^0(\bar{\Omega})$ to the function ${\cal I}_Hv = \sum_{i=1}^{n_I}v_i\hat\varphi_i\in U_H$, uniquely identified by the vector $\mathbf v=\{v_i\}_{i=1}^{n_I}$, where $v_i=v(\mathbf x_i)$.

In order to evaluate the function ${\cal I}_HBu^{\cal N\!N}$ at the required quadrature points $\{\mathbf x_j^q, j=1,...,n_q\}$ during the loss function computation or the final evaluation, one just needs to compute the quantities
\[
{\cal I}_HBu^{\cal N\!N}(\mathbf x_j^q) = \sum_{i=1}^{n_I}\left(Bu^{\cal N\!N}\right)(\mathbf x_i)\hat\varphi_i(\mathbf x_j^q).
\]
In practice, it is more convenient to introduce the sparse matrix $M\in \mathbb R^{n_q\times n_I}$ such that $M_{i,j}=\hat\varphi_j(\mathbf x_i^q)$. This allows us to evaluate the interpolated function at each quadrature point with a matrix-vector multiplication as follows:
\begin{equation}\label{eq:matrix_interpolation}
\begin{bmatrix}
{\cal I}_HBu^{\cal N\!N}\left(\mathbf x_1^q\right) \\
{\cal I}_HBu^{\cal N\!N}\left(\mathbf x_2^q\right) \\
\vdots \\
{\cal I}_HBu^{\cal N\!N}\left(\mathbf x_{n_q}^q\right)
\end{bmatrix}
 = M 
 \begin{bmatrix}
\left(Bu^{\cal N\!N}\right)(\mathbf x_1) \\
\left(Bu^{\cal N\!N}\right)(\mathbf x_2) \\
\vdots \\
\left(Bu^{\cal N\!N}\right)(\mathbf x_{n_I})
\end{bmatrix}.
\end{equation}
In the same way, the derivatives of ${\cal I}_HBu^{\cal N\!N}$ can be computed, at the same points, as
\begin{equation}\label{eq:interp_deriv}
\frac{\partial^{\vert\alpha\vert} {\cal I}_HBu^{\cal N\!N}}{\partial x^\alpha}
(\mathbf x_j^q) = \sum_{i=1}^{n_I}\left(Bu^{\cal N\!N}\right)(\mathbf x_i)
\frac{\partial^{\vert\alpha\vert} \hat\varphi_i(\mathbf x_j^q)}{\partial x^\alpha},
\end{equation}
where $\alpha=(\alpha_1,...,\alpha_{n})\in \mathbb Z_+^{n}$ and $\vert\alpha\vert = \sum_{i=1}^{n}\alpha_i$. Defining the matrix $M^\alpha\in \mathbb R^{n_q\times n_I}$ such that $M_{i,j}^\alpha=\frac{\partial^{\vert\alpha\vert} \hat\varphi_j(\mathbf x_i^q)}{\partial x^\alpha}$, it is possible to compute all the required derivatives simply by replacing $M$ by $M^\alpha$ on the right hand side of  \eqref{eq:matrix_interpolation}. In this way, the VPINN derivatives can be computed without relying on automatic differentiation, further improving the method efficiency during training.

\section{An example of `spurious' neural network}\label{sect:stability_large_dataset}

Consider again the boundary-value problem \eqref{eq:sol0}, which admits the null solution.
We are interested in solving this problem by a plain PINN (VPINN) solver. To train the network, we choose a set of $n_{\mathcal S}$ control (quadrature) points, 
\[
\mathcal S = \{x_i: i=1,...,n_{\mathcal S}\} 
\]
satisfying $0\le x_1 < x_2 < ... < x_{n_{\mathcal S}-1} < x_{n_{\mathcal S}}\le 1$.
Using the architecture defined in \eqref{nn_main_formula}, it is possible to construct a ReLU neural network $w$ with just a single hidden layer and 3 neurons with the following weights:
\[
\mathbf{A}_1^j = \begin{bmatrix}
    1/h_j \\
    1/h_j \\
    1/h_j
\end{bmatrix},
\hspace{0.5cm}
\mathbf{b}_1^j =  \begin{bmatrix}
    (- \overline x_j+h_j)/h_j \\
    -\overline x_j /h_j \\
    (- \overline x_j-h_j)/h_j
\end{bmatrix},
\hspace{0.5cm}
\mathbf{A}_2^j = \begin{bmatrix}
    \frac1h_j & -\frac2h_j & \frac1h_j
\end{bmatrix},
\hspace{0.5cm}
\mathbf{b}_2^j = 0,
\]
where, for any fixed index $j\in\{1,...,n_{\mathcal S}-1\}$, we denote by $\overline x_j=\frac{x_j+x_{j+1}}2$ the mean of two consecutive nodes,  and by $h_j$ the difference between $\overline x_j$ and $x_j+\varepsilon_j$, for some $\varepsilon_j\in(0,\overline x_j-x_j)$.
The function represented by this set of weights is:
\begin{equation}\label{eq:counter_definition}
\begin{aligned}
w_j(x)= & \frac 1h_j \max\left(0, \frac xh_j + \frac{-\overline x_j+h_j}{h_j}\right) 
-\frac 2h_j \max\left(0, \frac xh_j + \frac{-\overline x_j}{h_j}\right) + \\
&+\frac 1h_j \max\left(0, \frac xh_j + \frac{-\overline x_j-h_j}{h_j}\right)\\
=&\left\{\begin{aligned}
&-\frac{\overline x_j-h_j}{h_j^2}+\frac{x}{h_j^2} \hspace{0.5cm} && \text{ if } x\in(\overline x_j-h_j, \overline x_j] \,,\\
& \frac{\overline x_j+h_j}{h_j^2}-\frac{x}{h_j^2} \hspace{0.5cm} && \text{ if } x\in(\overline x_j, \overline x_j + h_j) \,, \\
& 0 &&   \text{ otherwise}.
\end{aligned}\right.
\end{aligned}
\end{equation}
An example of such a function is shown in Fig. \ref{fig:infsup_counterexample}.
\begin{figure}[t!]
\centering 
\captionsetup{justification=centering}
  \includegraphics[width=0.6\linewidth]{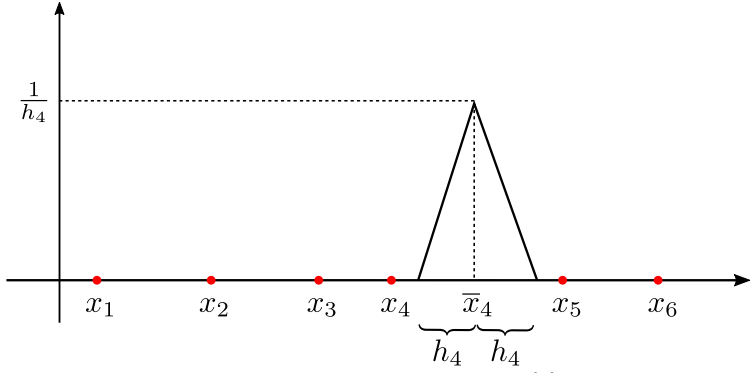} 
  \caption{Graphical representation of a function $u_4^{\cal N\!N}$ defined in \eqref{eq:counter_definition}.}
  \label{fig:infsup_counterexample}
\end{figure}

It is easily seen that $w_j(x_i)=w_j'(x_i)=w_j''(x_i)=0$ for any $i=1,...,n_{\mathcal S}$, therefore the PINN (VPINN) loss function is exactly equal to 0. However, this does not ensure the accuracy of the approximation, in fact $\Vert w_j\Vert_{L^1(\Omega)}=\Vert w_j-u\Vert_{L^1(\Omega)}=1$.

Note that it is possible to define larger networks with analogous properties. Moreover, the same phenomenon can be observed with any sigmoid activation function, exploiting the fact that the accuracy in the evaluation of both the activation function and the loss function is bounded by machine precision.
We also highlight that the phenomenon can be partially alleviated by introducing a regularization term in the loss function; however, it adds noise to the optimization process, possibly resulting in losses of accuracy when the PDE solution is characterized by large gradients.

This proves that it is not possible to guarantee the inf-sup stability for standard PINNs or VPINNs, and that spurious modes cannot be controlled simply minimizing the loss function. On the contrary, the interpolation operator proposed in this paper acts as a VPINN stabilizer, preventing the onset of spurious components in the solution.

}